\documentclass[10pt]{article}

\usepackage[]{amsmath, amssymb, theorem}
\usepackage{amsfonts}

\usepackage{theorem}
{\theorembodyfont{} \newtheorem{proposition}{Proposition}[section]}
{\theorembodyfont{} \newtheorem{theorem}[proposition]{Theorem}}
{\theorembodyfont{} \newtheorem{lemma}[proposition]{Lemma}}
{\theorembodyfont{\rm} \newtheorem{definition}[proposition]{Definition}}
{\theorembodyfont{} \newtheorem{corollary}[proposition]{Corollary}}
{\theorembodyfont{\rm} }
{\theorembodyfont{\rm} \newtheorem{remark}[proposition]{Remark}}
{\theorembodyfont{\rm} }
{\theorembodyfont{} }
{\theorembodyfont{} }
{\theorembodyfont{} }
{\theorembodyfont{} }

\newcommand\qed{\hfill $\Box$}
\newcommand{\QCBZ}{\mathsf{QCB_0}}
\newcommand{\CBZ}{\mathsf{CB_0}}

\newcommand{\Pomega}{{P\hspace*{-1pt}\omega}}
\newcommand{\bfSig}{\mathbf{\Sigma}}
\newcommand{\bfPi}{\mathbf{\Pi}}
\newcommand{\bfGamma}{\mathbf{\Gamma}}
\newcommand{\dom}{\mathit{dom}}
\newcommand{\rng}{\mathit{rng}}

\newcommand{\graph}{\mathit{graph}}

\newcommand{\calN}{\mathcal{N}}
\newcommand{\calO}{\mathcal{O}}

\newcommand{\IQ}{\mathbb{Q}}
\newcommand{\IS}{\mathbb{S}}
\newcommand{\IR}{\mathbb{R}}
\title{Towards a Descriptive Theory of cb$_0$-Spaces}

\author{Victor Selivanov\\A.P. Ershov
Institute of Informatics
Systems SB RAS\\
Novosibirsk, Russia}

\date{\today}

\begin{document}

\maketitle

\begin{abstract}
The paper  tries to extend results of the classical Descriptive Set
Theory to as many  countably based $T_0$-spaces (cb$_0$-spaces) as
possible. Along with extending some central facts about Borel, Luzin
and Hausdorff hierarchies of sets we consider also the more general
case of $k$-partitions. In particular, we investigate the difference
hierarchy of $k$-partitions and the fine hierarchy closely related
to the Wadge hierarchy.

\end{abstract}

%%%%%%%%%%%%%%%%%%%%%%%%%%%%%%%%%%%%%%%%%%%%%%%%%%%%%%%%%%%%%%
%
%
%%%%%%%%%%%%%%%%%%%%%%%%%%%%%%%%%%%%%%%%%%%%%%%%%%%%%%%%%%%%%%

\section{Introduction}\label{in}

Classical Descriptive  Set Theory  \cite{ke95} is an important field
of mathematics with numerous applications. It investigates
descriptive complexity of sets, functions and equivalence relations
in Polish (i.e., separable complete metrizable) spaces.

Although  Polish spaces are sufficient for many fields of classical
mathematics, they are certainly not sufficient for many fields of
Theoretical Computer Science where non-Hausdorff spaces (in
particular, $\omega$-continuous domains) and non-countably-based
spaces (in particular, Kleene-Kreisel continuous functionals) are of
central importance. For this reason, the task of extending classical
Descriptive Set Theory (DST) to as many non-Polish spaces as
possible attracted attention of several researchers.

Some  parts of DST for $\omega$-continuous domains (which are
typically non-Polish) were developed in \cite{s04,s05a,s06,s08}. In
\cite{br} a good deal of DST was developed for the so called
quasi-Polish spaces (see the next section for a  definition of this
class of cb$_0$-spaces) which include both the Polish spaces and the
$\omega$-continuous domains. For some attempts to develop DST for
non-countably based spaces see e.g. \cite{jr82,ms10,fhk11,p12,pb13}.

In this paper,  we try to develop  DST for some classes of
cb$_0$-spaces  beyond the class of quasi-Polish spaces. As is usual
in classical DST, we put emphasis on the ``ininitary version'' of
hierarchy theory where people are concerned with transfinite (along
with finite) levels of hierarchies. The ``finitary'' version where
people concentrate on the finite levels of hierarchies has a special
flavor and is relevant to several fields of Logic and Computation
Theory; it was systematized in \cite{s06,s08a,s12}.

We extend some well  known facts about classical hierarchies in
Polish spaces to natural classes of cb$_0$-spaces. Namely, we show
that some levels of hierarchies of cb$_0$-spaces introduced in
\cite{scs13,scs14} provide natural examples of classes of
cb$_0$-spaces with reasonable DST (in particular, the classical
Suslin, Hausdorff-Kuratowski and non-collapse theorems for the
Borel, Luzin and Hausdorff hierarchies are true for such spaces).
This portion of our results are technically easy and follow rather
straightforwardly from the classical DST and some notions and
results in \cite{mss12,scs13}.

Along with  the classical hierarchies of sets we are interested also
in the difference and fine hierarchies of $k$-partitions
\cite{s06,s07,s07a,s08,s08a,s11} which seem to be natural,
non-trivial and useful generalization of the corresponding
hierarchies of sets. Also, along with the classical Wadge
reducibility \cite{wad72,wad84,vw76} we discuss its extension to
$k$-partitions \cite{he93,he96,s07a}, and  some their weaker
versions introduced and studied in \cite{am03,an06,mr09,mss12}.

Already the  extension of the Hausdorff difference hierarchy to
$k$-partitions is a non-trivial task. The general ``right'' finitary
version of this hierarchy was found only recently in \cite{s12},
although for some particular cases it was already in our previous
publications. The general ``right'' infinitary version of this
hierarchy is new here, although the definition adequate for bases
with the $\omega$-reduction property was also found earlier
\cite{s07,s07a,s08}. That the definition in this paper is right
follows from the nice properties of the difference hierarchy of
$k$-partition (in particular, the natural version of the
Hausdorff-Kuratowski theorem).

The situation with the fine hierarchy (which aims to extend the
Wadge hierarchy to the case of sets and $k$-partitions in arbitrary
spaces) is even more complicated. This task is not obvious even for
the case of sets because the Wadge hierarchy is developed so far
only for the Baire space (and some of its close relatives) in terms
of $m$-reducibility by continuous functions and with a heavy use of
Martin determinacy theorem \cite{wad72,wad84,vw76}. As a result,
there is no clear explicit description of levels of the hierarchy in
terms of set-theoretic operations which one could try to extend to
other spaces (more precisely, some rather indirect descriptions
presented in \cite{wad84,lo83} strongly depend on the
$\omega$-reduction property of the open sets which usually fails in
non-zero-dimensional spaces). Probably, that was the reason why some
authors tried to obtain alternative characterizations of levels of
the Wadge hierarchy \cite{lo83,du01}. In a series of our papers (see
e.g. \cite{s06,s08a}) a characterization of levels of an abstract
version of Wadge hierarchy in the finitary case was achieved that
was extended in \cite{s12} to the case of $k$-partitions. Here we
develop an infinitary version of this approach and try to explain
why the corresponding hierarchy is the ``right'' extension of the
Wadge hierarchy to arbitrary spaces and to the $k$-partitions. Since
the notation and full proofs in this context are extremely involved,
we concentrate here only on the partitions of finite Borel rank and
avoid the complete proofs of some complicated results, giving only
precise formulations and hints of proofs with references to closely
related earlier proofs in the finitary context.

After  recalling some notions and known facts in the next section,
we discuss some basic properties of Borel and Luzin hierarchies in
cb$_0$-spaces in Section \ref{sub:embed}. In Section \ref{opens} we
establish some basic  facts on the difference hierarchies of
$k$-partitions in  cb$_0$-spaces. The main result here is the
Hausdorff-Kuratowski theorem for $k$-partitions in quasi-Polish
spaces. In Section \ref{fhs} we extend the difference hierarchies of
$k$-partitions to the the fine hierarchies of $k$-partitions. In
particular, we extend the Hausdorff-Kuratowski theorem to the fine
hierarchy. We conclude in Section~\ref{con} with sketching a
possible further research on extending the classical Descriptive Set
Theory.

%%%%%%%%%%%%%%%%%%%%%%%%%%%%%%%%%%%%%%%%%%%%%%%%%%%%%%%%%%%%%%%%
%
%
%%%%%%%%%%%%%%%%%%%%%%%%%%%%%%%%%%%%%%%%%%%%%%%%%%%%%%%%%%%%%%%%

\section{Notation and preliminaries}\label{prelim}

In this section we recall some notation, notions and results used in
the subsequent sections.

\subsection{cb$_0$-Spaces and qcb$_0$-spaces}\label{spaces}

Here we recall some topological notions and facts relevant to this
paper.

We freely use the standard set-theoretic notation like
$\dom(f),\rng(f)$ and $\graph(f)$ for the domain, range and graph of
a function $f$, respectively, $X\times Y$ for the Cartesian product,
and $P(X)$ for the set of all
subsets of $X$. For $A\subseteq X$, $\overline{A}$ denotes the
complement $X\setminus A$ of $A$ in $X$. We identify the set of
natural numbers with the first infinite ordinal $\omega$. The first
uncountable ordinal is denoted by $\omega_1$. The notation $f:X\to
Y$ means that $f$ is a total function from a set $X$ to a set $Y$.

We assume the reader to be familiar with the basic notions of
topology \cite{en89}. The collection of all open subsets of a
topological space $X$ (i.e.\ the topology of $X$) is denoted by
$\calO(X)$; for the underlying set of $X$ we will write $X$ in abuse
of notation. We will often abbreviate ``topological space'' to
``space''. A space is \emph{zero-dimensional} if it has a basis of
clopen sets. Recall that a \emph{basis} for the topology on $X$ is a
set $\cal B$ of open subsets of $X$ such that for every $x\in X$ and
open $U$ containing $x$ there is $B\in \cal B$ satisfying $x\in
B\subseteq U$.

Let $\omega$ be the space of non-negative integers with the
discrete topology. Of course, the spaces
$\omega\times\omega=\omega^2$, and $\omega\sqcup\omega$ are
homeomorphic to $\omega$, the first homeomorphism is realized by
the Cantor pairing function $\langle \cdot,\cdot\rangle$.

Let $\calN=\omega^\omega$ be the set of all infinite
sequences of natural numbers (i.e., of all functions $\xi \colon
\omega \to \omega$). Let $\omega^*$ be the set of finite sequences
of elements of $\omega$, including the empty sequence. For
$\sigma\in\omega^*$ and $\xi\in\calN$, we write
$\sigma\sqsubseteq \xi$ to denote that $\sigma$ is an initial
segment of the sequence $\xi$. By $\sigma\xi=\sigma\cdot\xi$ we
denote the concatenation of $\sigma$ and $\xi$, and by
$\sigma\cdot\calN$ the set of all extensions of $\sigma$ in
$\calN$. For $x\in\calN$, we can write
$x=x(0)x(1)\dotsc$ where $x(i)\in\omega$ for each $i<\omega$. For
$x\in\calN$ and $n<\omega$, let $x^{<n}=x(0)\dotsc x(n-1)$
denote the initial segment of $x$ of length $n$. Notations in the
style of regular expressions like $0^\omega$, $0^\ast 1$ or
$0^m1^n$ have the obvious standard meaning.

By endowing $\calN$ with the product of the discrete
topologies on $\omega$, we obtain the so-called \emph{Baire space}.
The product topology coincides with the topology
generated by the collection of sets of the form
$\sigma\cdot\calN$ for $\sigma\in\omega^*$. The Baire space
is of primary importance for Descriptive Set Theory and Computable Analysis.
The importance stems from
the fact that many countable objects are coded straightforwardly
by elements of $\calN$, and it has very specific topological
properties. In particular, it is a perfect zero-dimensional space
and the spaces $\calN^2$, $\calN^\omega$,
$\omega\times\calN=\calN\sqcup{\mathcal
N}\sqcup\cdots$ (endowed with the product topology) are all
homeomorphic to $\calN$. Let $(x,y)\mapsto\langle
x,y\rangle$ be a homeomorphism between $\calN^2$ and
$\calN$.
The subspace $\mathcal{C}:=2^\omega$ of $\calN$ formed by
the infinite binary strings (endowed with the relative topology
inherited from $\calN$) is known as the \emph{Cantor space}.

The  \emph{Sierpinski space} $\mathbb{S}$ is the two-point set
$\{\bot,\top\}$ where the set $\{\top\}$ is open but not closed. The
space $\Pomega$ is formed by the set of subsets of $\omega$ equipped
with the Scott topology. A countable base of the Scott topology is
formed by the sets $\{A\subseteq\omega\mid F\subseteq A\}$, where
$F$ ranges over the finite subsets of $\omega$. Note that
$\Pomega=\calO(\omega)$. As is well-known \cite{g03}, $\Pomega$ is
universal for cb$_0$-spaces:

\begin{proposition}\label{inj2}
 A topological space $X$ embeds into $\Pomega$ iff $X$ is a $cb_0$-space.
\end{proposition}

Remember  that a space $X$ is \emph{Polish} if it is countably based
and metrizable with a metric $d$ such that $(X,d)$ is a complete
metric space. Important examples of Polish spaces are $\omega$,
$\calN$, $\mathcal{C}$, the space of reals $\IR$ and its Cartesian
powers $\IR^n $ ($ n < \omega $), the closed unit interval $ [0,1]
$, the Hilbert cube $ [0,1]^\omega $ and the Hilbert space
$\IR^\omega $. Simple examples of non-Polish spaces are $\IS$,
$\Pomega$ and the space $\IQ$ of rationals.

A space $X$ is {\em quasi-Polish} \cite{br} if it is
countably based and quasi-metrizable with a quasi-metric $d$ such
that $(X,d)$ is a complete quasi-metric space. A {\em
quasi-metric} on $X$ is a function from $X\times X$ to the
nonnegative reals such that $d(x,y)=d(y,x)=0$ iff $x=y$, and
$d(x,y)\leq d(x,z)+d(z,y)$. Since for the quasi-metric spaces
different notions of completeness and of a Cauchy sequence are
considered, the definition of quasi-Polish spaces should be made
more precise (see \cite{br} for additional details).  We skip
these details because we will  in fact use another characterization
of these spaces given below.
Note that the spaces $\IS$, $\Pomega$ are quasi-Polish while the space $\IQ$ is not.

A \emph{representation} of a space $X$ is a surjection of a subspace
of the Baire space $\calN$ onto $X$. A basic notion of Computable
Analysis is the notion of admissible representation. A
representation $\delta$ of $X$ is \emph{admissible}, if it is
continuous and any continuous function $\nu:Z \to X$ from a subset
$Z\subseteq\calN$ to $X$ is continuously reducible to $\delta$, i.e.
$\nu=\delta\circ g$ for some continuous function $g:Z \to \calN$. A
topological space is \emph{admissibly representable}, if it has an
admissible representation.

The  notion of admissibility was introduced in \cite{kw85} for
representations of cb$_0$-spaces (in a different but equivalent
formulation) and was extensively studied by many authors. In
\cite{sch:ext,sch:phd} the notion was extended to non-countably
based spaces and a nice characterization of the admissibly
represented spaces was achieved. Namely, the admissibly represented
sequential topological spaces coincide with the qcb$_0$-spaces,
i.e., $T_0$-spaces which are topological quotients of countably
based spaces.

In \cite{br} the following important characterization of
quasi-Polish spaces in terms of Borel hierarchy was obtained.

\begin{proposition}\label{pi2}
A space is quasi-Polish iff it is homeomorphic to a ${\bf
\Pi}^0_2$-subset of $P\omega$ with the induced topology.
\end{proposition}

%%%%%%%%%%%%%%%%%%%%%%%%%%%%%%%%%%%%%%%%%%%%%%%%%%%%%%%%%%%%%%%%%
%

\subsection{Hierarchies of sets}\label{fam}

Here we briefly recall definitions and some properties of Borel and
Luzin hierarchies in arbitrary topological spaces.

A {\em pointclass} in a space $ X $ is simply a collection $
\mathbf{\Gamma}(X) $ of subsets of $ X $. A {\em family of
pointclasses} \cite{s13} is a family $
\mathbf{\Gamma}=\{\mathbf{\Gamma}(X)\} $ indexed by arbitrary
topological spaces $X$ such that each $ \mathbf{\Gamma}(X) $ is a
pointclass on $ X $ and $ \mathbf{\Gamma} $ is closed under
continuous preimages, i.e. $ f^{-1}(A)\in\mathbf{\Gamma}(X) $ for
every $ A\in\mathbf{\Gamma}(Y) $ and every continuous function $ f
\colon X\to Y $. A basic example of a family of pointclasses is
given by the family $\mathcal{O}=\{\mathcal{O}(X)\}$ of the
topologies of all the spaces $X$.

We will use some operations on families of pointclasses. First, the usual
set-theoretic operations will be applied to the families of
pointclasses pointwise: for example, the union $\bigcup_i
\mathbf{\Gamma}_i$ of the families of pointclasses
$\mathbf{\Gamma}_0,\mathbf{\Gamma}_1,\ldots$ is defined by
$(\bigcup_i\mathbf{\Gamma}_i)(X)=\bigcup_i\mathbf{\Gamma}_i(X)$.

Second, a large class of such operations is induced by the
set-theoretic operations of L.V. Kantorovich and E.M. Livenson (see
e.g. \cite{s13} for the general definition). Among them are the
operation $\bfGamma\mapsto\bfGamma_\sigma$, where
$\mathbf{\Gamma}(X)_\sigma$ is the set of all countable unions of
sets in $\mathbf{\Gamma}(X)$, the operation
$\mathbf{\Gamma}\mapsto\mathbf{\Gamma}_\delta$, where
$\mathbf{\Gamma}(X)_\delta$ is the set of all countable
intersections of sets in $\mathbf{\Gamma}(X)$, the operation
$\mathbf{\Gamma}\mapsto\mathbf{\Gamma}_c=\check{\mathbf{\Gamma}}$,
where $\mathbf{\Gamma}(X)_c$ is the set of all complements of sets
in $\mathbf{\Gamma}(X)$, the operation
$\mathbf{\Gamma}\mapsto\mathbf{\Gamma}_d$, where $\bfGamma(X)_d$ is
the set of all differences of sets in $\bfGamma(X)$, the operation
$\bfGamma\mapsto\bfGamma_\exists$ defined by
$\bfGamma_\exists(X):=\{\exists^\mathcal{N}(A)\mid
A\in\bfGamma(\calN\times X)\}$, where
$\exists^\mathcal{N}(A):=\{x\in X\mid \exists p\in\calN.(p,x)\in
A\}$ is the projection of $A\subseteq\calN\times X$ along the axis
$\calN$, and finally the operation $\bfGamma\mapsto\bfGamma_\forall$
defined by $\bfGamma_\forall(X):=\{ \forall^\mathcal{N}(A) \mid
A\in\bfGamma(\calN\times X)\}$, where
$\forall^\mathcal{N}(A):=\{x\in X\mid \forall p\in\calN.(p,x)\in
A\}$.

The  operations on families of pointclasses enable to provide short
uniform descriptions of the classical hierarchies in arbitrary
spaces. E.g., the Borel  hierarchy is the family of pointclasses
$\{\bfSig^0_\alpha\}_{ \alpha<\omega_1}$ defined by induction on
$\alpha$ as follows \cite{s06,br}: $\bfSig^0_0(X):=\{\emptyset\}$,
$\bfSig^0_1 := \mathcal{O}$, $\bfSig^0_2 := (\bfSig^0_1)_{d\sigma}$,
and  $\bfSig^0_\alpha(X)
:=(\bigcup_{\beta<\alpha}\bfSig^0_\beta(X))_{c\sigma}$ for
$\alpha>2$. The sequence $\{\bfSig^0_\alpha(X)\}_{ \alpha<\omega_1}$
is called \emph{the Borel hierarchy} in $X$. We also let
$\bfPi^0_\beta(X) := (\bfSig^0_\beta(X))_c $ and
$\mathbf{\Delta}^0_\alpha(X) := \bfSig^0_\alpha(X) \cap
\bfPi^0_\alpha (X)$. The classes
$\bfSig^0_\alpha(X),\bfPi^0_\alpha(X),{\bf\Delta}^0_\alpha(X)$ are
called the \emph{levels} of the Borel hierarchy in $X$.

We  recall an important structural property of
$\mathbf{\Sigma}$-levels of the Borel hierarchy. Let $\Gamma$ be a
family of pointclasses. A pointclass $\Gamma(X)$ has the {\em
$\omega$-reduction property} if  for each countable sequence $A_0,
A_1,\ldots$ in $\Gamma(X)$ there is a countable sequence $D_0,
D_1,\ldots$ in $\Gamma(X)$ such that $D_i\subseteq A_i$, $D_i\cap
D_j=\emptyset$ for all $i\not=j$ and
$\bigcup_{i<\omega}D_i=\bigcup_{i<\omega}A_i$.

\begin{proposition}\label{uniform}
For any space $X$ and any $2\leq\alpha<\omega_1$,
${\bf\Sigma}^0_\alpha(X)$  has the $\omega$-reduction  properties.
If $X$ is zero-dimensional, the same holds for the class
${\bf\Sigma}^0_1(X)$ of open sets.
\end{proposition}

The {\em hyperprojective hierarchy} is the family of pointclasses
 $\{\bfSig^1_\alpha\}_{ \alpha<\omega_1}$ defined by induction on $\alpha$ as follows:
 $\bfSig^1_0=\bfSig^0_2$,
 $\bfSig^1_{\alpha+1}=(\bfSig^1_\alpha)_{c\exists}$,
 $\bfSig^1_\lambda=(\bfSig^1_{<\lambda})_{\delta \exists}$,
 where $\alpha,\lambda<\omega_1$, $\lambda$ is a limit ordinal,
 and $\bfSig^1_{<\lambda}(X):=\bigcup_{\alpha<\lambda}\bfSig^1_\alpha(X)$.

In this way, we obtain for any topological space $X$ the sequence
$\{\bfSig^1_\alpha(X)\}_{\alpha<\omega_1}$, which we call here
\emph{the hyperprojective hierarchy in $X$}. The pointclasses
$\bfSig^1_\alpha(X)$, $\bfPi^1_\alpha(X):=(\bfSig^1_\alpha(X))_c$
and
$\mathbf{\Delta}^1_\alpha(X):=\bfSig^1_\alpha(X)\cap\bfPi^1_\alpha(X)$
are called \emph{levels of the hyperprojective hierarchy in $X$}.
The finite non-zero levels of the hyperprojective hierarchy coincide
with the corresponding levels of the Luzin's projective hierarchy
\cite{br,scs13}. The class of  \emph{hyperprojective sets} in $X$ is
defined as the union of all levels of the hyperprojective hierarchy
in $X$. For more information on the hyperprojective hierarchy see
\cite{ke83,ke95,scs14}. Below we will also consider some other
hierarchies, in particular the Hausdorff difference hierarchy.

\subsection{$k$-Partitions and hierarchies over well posets}

Here we discuss a more general notion of a hierarchy (compared with
the notion of hierarchy of sets \cite{s08a,s12}) which applies, in
particular, to the hierarchies of $k$-partitions.

Let  $2\leq k<\omega$. By a {\em $k$-partition of a space $X$} we
mean a function $A:X\to k=\{0,\ldots,k-1\}$ often identified with
the sequence $(A_0,\ldots,A_{k-1})$ where $A_i=A^{-1}(i)$.
Obviously, 2-partitions of $X$ are identified with the subsets of
$X$ using the characteristic functions. The set of all
$k$-partitions of $X$ is denoted $k^X$. For
$\mathbf{\Gamma}\subseteq P(X)$, let $(\mathbf{\Gamma})_k$ denote
the set of $k$-partitions $A\in k^X$ such that
$A_0,\ldots,A_{k-1}\in\mathbf{\Gamma}(X)$. In particular,
$(\mathbf{\Sigma}_{<\omega}(X))_k$ is the set of $k$-partitions of
finite Borel rank which will be considered in Section \ref{fhpart}.

The Wadge  reducibility on subsets of $X$ is naturally extended on
$k$-partitions: for $A,B\in k^X$, $A\leq_W B$ means that $A=B\circ
f$ for some continuous function $f$ on $X$, In this way, we obtain
the preorder $(k^X;\leq_W)$ which for $k\geq3$ turns out much more
complicated than the structure of Wadge degrees, even for simple
case $X=\mathcal{N}$ \cite{he93,he96,s06,s11}.

To find the  ``right'' extensions of the classical difference and
Wadge hierarchies from the case of sets to the case of
$k$-partitions is a quite non-trivial task. A reason is that levels
of hierarchies of sets are always semi-well-ordered by inclusion (in
particular, there are no three levels which are pairwise
incomparable by inclusion) while the structure of hierarchies of
$k$-partitions for $k\geq3$ is usually more complicated than the
structure of the hierarchies of sets (in particular, for $k\geq3$
the poset of levels of difference hierarchies of $k$-partitions
under inclusion usually has antichains with any finite number of
elements).

Here we recall from \cite{s12} a very general notion of a hierarchy
that covers all hierarchies  we discuss in this paper.

 \begin{definition}\label{hierpart0}
\begin{enumerate}\itemsep-1mm
 \item For any poset $P$ and any set $A$, by a $P$-hierarchy in
$A$ we mean a family $\{H_p\}_{p\in P}$ of subsets of $A$ such
that $p\leq q$ implies $H_p\subseteq H_q$.
  \item Levels (resp. constituents) of a $P$-hierarchy
  $\{H_p\}$ are the sets $H_{p_0}\cap\cdots\cap
H_{p_n}$ (resp. the sets $C_{p_0,\ldots,p_n}=(H_{p_0}\cap\cdots\cap
H_{p_n})\setminus\bigcup\{H_q\mid q\in
P\setminus\uparrow\{p_0,\ldots,p_n\}\}$)
 where $n\geq0$ and $\{p_0,\ldots,p_n\}$ is an antichain in
$P$.
 \item A $P$-hierarchy $\{H_p\}$ is precise if $p\leq q$ is
equivalent to $H_p\subseteq H_q$.
  \end{enumerate}
 \end{definition}

Note that the classical hierarchies of sets are obtained from the
above definition if $A=P(X)$ and $P=\bar{2}\cdot\eta$ is the poset
obtained by replacing any element of the ordinal $\eta$ by an
antichain with two elements, and that the notion of preciseness
extends the non-collapse property of hierarchies. Note that levels
of the classical hierarchies coincide with levels in the sense of
the definition above. The constituents of say, Borel hierarchy, are
$\mathbf{\Sigma}^0_\alpha\setminus\mathbf{\Pi}^0_\alpha$,
$\mathbf{\Pi}^0_\alpha\setminus\mathbf{\Sigma}^0_\alpha$,
$\mathbf{\Delta}^0_{\alpha+1}
\setminus(\mathbf{\Sigma}^0_\alpha\cup\mathbf{\Pi}^0_\alpha)$, and
$\mathbf{\Delta}^0_{\lambda}
\setminus\bigcup_{\alpha<\lambda}\mathbf{\Sigma}^0_\alpha$, where
$\lambda$ is a limit countable ordinal.

As it was already mentioned, for hierarchies of $k$-partitions
(obtained when $A=k^X$) we cannot hope to deal only with
semi-well-ordered posets $P=\bar{2}\cdot\eta$ in the definition
above. Fortunately, a slight weakening of this property is
sufficient for our purposes: we can confine ourselves with the so
called well posets (wpo) or, more generally well preorders (wqo).
Recall that a wqo is a preorder $P$ that has neither infinite
descending chains nor infinite antichains. The theory of wqo (widely
known as the wqo-theory) is a well developed field with several deep
results and applications, see e.g. \cite{kru72}. It is also of great
interest to hierarchy theory. An important role in wqo-theory
belongs to the rather technical notion of a better preorder (bqo).
Bqo's form a subclass of wqo's with good closure properties.

Note that if $P$ is a wpo then the structure $(\{H_p\mid p\in
P\};\subseteq)$ of levels of a $P$-hierarchy under inclusion is also
a wpo, hence some important features of the hierarchies of sets hold
also for the hierarchies of partitions. Moreover, for such
hierarchies we have some  important properties of constituents, in
particular the constituents form a partition of the set
$\bigcup\{H_p\mid p\in P\}$ (see also Section 7 of \cite{s12} for
additional details).

Although well posets are very simple compared with arbitrary posets,
they are much more complicated than the semi-well-orders which
essentially reduce to the ordinals. Obviously, there are a lot
isomorphism types of well posets of a fixed rank. Below we consider
some examples of well posets suitable for naming the levels of the
difference and fine hierarchies of $k$-partitions.

%%%%%%%%%%%%%%%%%%%%%%%%%%%%%%%%%%%%%%%%%%%%%%%%%%%%%%%

\subsection{Hierarchies of cb$_0$-spaces and qcb$_0$-spaces}\label{qcb}

Here we recall some classifications of qcb$_0$-spaces induced by the
classical hierarchies of sets.

For any representation $\delta$ of a space $X$, let
$\mathit{EQ}(\delta):=\{\langle p,q\rangle \in \calN \mid p,q\in
\dom(\delta)\wedge\delta(p)=\delta(q)\}$. Let $\mathbf{\Gamma}$ be a
family of pointclasses. A qcb$_0$-space $X$ is called
\emph{$\mathbf{\Gamma}$-representable}, if $X$ has an admissible
representation $\delta$ with
$\mathit{EQ}(\delta)\in\mathbf{\Gamma}(\calN)$. The class of all
$\mathbf{\Gamma}$-representable spaces is denoted
$\QCBZ(\mathbf{\Gamma})$. A cb$_0$-space $X$ is called a
\emph{$\mathbf{\Gamma}$-space}, if $X$ is homeomorphic to a
$\mathbf{\Gamma}$-subspace of $\Pomega$. The class of all
$\mathbf{\Gamma}$-spaces is denoted $\CBZ(\mathbf{\Gamma})$.

These  notions from \cite{scs13} enable to transfer hierarchies of
sets to the corresponding hierarchies of qcb$_0$-spaces. In
particular, we arrive at the following definition.

\begin{definition}\label{def:hierqcb}
 The sequence
$\{\CBZ(\bfSig^0_\alpha)\}_{\alpha<\omega_1}$ (resp. the sequence
$\{\QCBZ(\bfSig^0_\alpha)\}_{\alpha<\omega_1}$) is called the
\emph{Borel hierarchy} of cb$_0$-spaces (resp. of qcb$_0$-spaces).
By {\em levels} of this hierarchy we mean the classes
$\CBZ(\bfSig^0_\alpha)$ as well as the classes
$\CBZ(\bfPi^0_\alpha)$ and $\CBZ(\mathbf{\Delta}^0_\alpha)$. In a
similar way one can define the hyperprojective hierarchies of
cb$_0$- and of qcb$_0$-spaces.
\end{definition}

The following fact from \cite{scs13,scs14} shows that the introduced
hierarchies agree on  cb$_0$-spaces:

\begin{proposition}\label{p:equiv:Gamma-representable}
For any $\mathbf{\Gamma} \in \{\bfPi^0_2,\bfSig^0_\beta,
\bfPi^0_\beta, \bfSig^1_\alpha, \bfPi^1_\alpha \mid
1\leq\alpha<\omega_1, 3\leq\beta<\omega_1\}$, we have
$\QCBZ(\mathbf{\Gamma}) \cap \CBZ = \CBZ(\mathbf{\Gamma})$, where
$\CBZ$ is the class of all cb$_0$-spaces.
\end{proposition}

Note that, by Proposition \ref{pi2},  $\CBZ(\bfPi^0_2)$ coincides
with the class of quasi-Polish spaces.

%%%%%%%%%%%%%%%%%%%%%%%%%%%%%%%%%%%%%%%%%%%%%%%%%%%%%%%%%%%%%%%%%

\section{Borel and Luzin hierarchies}\label{sub:embed}

In this section we extend some classical facts on the Borel and
Luzin hierarchies in Polish spaces on larger classes of
cb$_0$-spaces.

\subsection{Some reducibilities and isomorphisms}

Here we provide some information on versions of the Wadge
reducibility and of the notion of homeomorphisms relevant to this
paper.

Let \( \mathbf{\Gamma} \) be a family of pointclasses and $X,Y$ be
topological spaces. By \( \mathbf{\Gamma}(X,Y) \) we denote the
class of functions $f:X\to Y$ such that
$f^{-1}(A)\in\mathbf{\Gamma}(X)$ whenever $A\in\mathbf{\Gamma}(Y)$.
A set $A\subseteq X$ is {\em \( \mathbf{\Gamma} \)-reducible} to a
set $B\subseteq X$ (in symbols, $A\leq_\mathbf{\Gamma}B$) if
$A=f^{-1}(B)$ for some \(f\in \mathbf{\Gamma}(X,X) \).

Note that the \( \mathbf{\Sigma}^0_1 \)-functions coincide with the
continuous functions and the \( \mathbf{\Sigma}^0_1 \)-reducibility
coincides with the classical Wadge reducibility. \(
\mathbf{\Sigma}^0_\alpha \)-Functions and \(
\mathbf{\Sigma}^0_\alpha \)-reducibilities were investigated in
\cite{am03,an06,mr09}.

We say that  topological spaces \( X , Y \) are \emph{\(
\mathbf{\Gamma} \)-isomorphic} if there is a bijection $f$ between
$X$ and $Y$ such that $f\in\mathbf{\Gamma}(X,Y)$ and
$f^{-1}\in\mathbf{\Gamma}(Y,X)$. It is a classical fact of
Descriptive Set Theory that  every two uncountable Polish spaces \(
X , Y \) are \( \mathbf{\Delta}^1_1 \)-isomorphic (see e.g.\
\cite[Theorem 15.6]{ke95}). The next result from \cite{mss12}
extends this fact to the context of uncountable quasi-Polish spaces
and computes an upper bound for the complexity of the
Borel-isomorphism.

\begin{proposition}\label{propgeneralhomeo}
Let  \( X,Y \) be two uncountable quasi-Polish spaces. Then  $X$ and
$Y$ are \( \mathbf{\Delta}^0_{<\omega} \)-isomorphic. If the
inductive dimensions $dim(X)$, $dim(Y)$ of  $X,Y$ are distinct from
$\infty$ then $X$ and $Y$ are \( \mathbf{\Delta}^0_3\)-isomorphic.
\end{proposition}

Let again $\mathbf{\Gamma}$ be a family of pointclasses. By a {\em
$\mathbf{\Gamma}$-family of pointclasses} we mean a family
$\{E(X)\}_X$ indexed by arbitrary spaces such that $E(X)$ is a
pointset in $X$, and $f^{-1}(A)\in E(X)$ for all $A\in E(Y)$ and
$f\in\mathbf{\Gamma}(X,Y)$. Obviously, the
$\mathbf{\Sigma}^0_1$-families of pointclasses are precisely the
``usual'' families of pointclasses.

\begin{lemma}\label{gfam}
Let $\mathbf{\Gamma}$ be a family of pointclasses. Then
$\mathbf{\Gamma}$ is a $\mathbf{\Gamma}$-family of pointclasses, any
continuous function $f:X\to Y$ is in $\mathbf{\Gamma}(X,Y)$, and any
$\mathbf{\Gamma}$-family of pointclasses is a family of
pointclasses.
\end{lemma}

{\bf Proof.}  The first assertion is obvious. For the second
assertion, let $f:X\to Y$ be continuous and
$A\in\mathbf{\Gamma}(Y)$. Since $\mathbf{\Gamma}$ is a family of
pointclasses, $f^{-1}(A)\in \mathbf{\Gamma}(X)$. Since $A$ was
arbitrary, $f\in\mathbf{\Gamma}(X,Y)$. The third assertion follows
from the second one.
 \qed

\begin{lemma}\label{gfam1}
Let $\alpha<\beta<\omega_1$.  Then any
$\mathbf{\Sigma}^0_\alpha$-function is a
$\mathbf{\Sigma}^0_\beta$-function, and any
$\mathbf{\Sigma}^1_\alpha$-function is a
$\mathbf{\Sigma}^1_\beta$-function.
\end{lemma}

{\bf Proof Hint.}  Proof is straightforward by induction on $\beta$,
so we consider only the first assertion for the case
$\beta=\alpha+1$, as an example. Let
$A\in\mathbf{\Sigma}^0_\beta(Y)$ and
$f\in\mathbf{\Sigma}^0_\alpha(X,Y)$. Then $A=\bigcup_n(Y\setminus
A_n)$ for some $A_0,A_1,\ldots\in\mathbf{\Sigma}^0_\alpha(Y)$. Then
$f^{-1}(A)=\bigcup_n(X\setminus
f^{-1}(A_n))\in\mathbf{\Sigma}^0_\beta(X)$, so $f$ is a
$\mathbf{\Sigma}^0_\beta$-function.
 \qed

\begin{proposition}\label{aliso}
\begin{enumerate}
 \item Let  $\mathbf{\Gamma}\in\{\mathbf{\Sigma}^0_\alpha,\mathbf{\Pi}^0_\alpha,
 \mathbf{\Sigma}^1_\beta,\mathbf{\Pi}^1_\beta\mid \omega\leq\alpha<\omega_1,1\leq\beta<\omega_1\}$
 and $X\in\CBZ(\mathbf{\Gamma})$. Then $X$ is $\mathbf{\Gamma}$-isomorphic
 to a subspace $S$ of $\mathcal{N}$ such that $S\in\mathbf{\Gamma}(\mathcal{N})$.
 \item Let  $\mathbf{\Gamma}\in\{\mathbf{\Sigma}^0_\alpha,\mathbf{\Pi}^0_\alpha,
 \mathbf{\Sigma}^1_\beta,\mathbf{\Pi}^1_\beta\mid 3\leq\alpha<\omega_1,1\leq\beta<\omega_1\}$
 and $X\in\CBZ(\mathbf{\Gamma})$, $dim(X)\not=\infty$.
 Then $X$ is $\mathbf{\Gamma}$-isomorphic to a subspace $S$ of $\mathcal{N}$
 such that $S\in\mathbf{\Gamma}(\mathcal{N})$.
 \end{enumerate}
\end{proposition}

{\bf Proof.}  Both items are checked in the same way, so consider
only the first one. Assume without loss of generality that
$X\in\mathbf{\Gamma}(P\omega)$. By Proposition
\ref{propgeneralhomeo}, the spaces  $P\omega$ and $\mathcal{N}$ are
\( \mathbf{\Delta}^0_{<\omega} \)-isomorphic, hence also \(
\mathbf{\Sigma}^0_{\omega} \)-isomorphic. By Lemma \ref{gfam1}
$P\omega$ and $\mathcal{N}$ are $\mathbf{\Gamma}$-isomorphic, let
$f:P\omega\to\mathcal{N}$ be a $\mathbf{\Gamma}$-isomorphism. Then
$f|_X$ is a $\mathbf{\Gamma}$-isomorphism between $X$ and
$S=f(X)\in\mathbf{\Gamma}(\mathcal{N})$.
 \qed

\subsection{Borel and Luzin hierarchies in cb$_0$-spaces}\label{propblh}

Here we  extend some well known facts on the Borel and Luzin
hierarchies in Polish and quasi-Polish spaces.

As is well-known,  any uncountable Polish (or quasi-Polish) space is
of continuum cardinality. The next fact extends this to many
cb$_0$-spaces:

\begin{proposition}\label{contin}
Any uncountable  space  $X$ in $\CBZ(\mathbf{\Sigma}^1_1)$ is of
continuum cardinality.
\end{proposition}

{\bf Proof.} By Proposition  \ref{aliso}, $X$ is
$\mathbf{\Sigma}^1_1$-isomorphic to a subspace $S$ of $\mathcal{N}$
such that $S\in\mathbf{\Sigma}^1_1(\mathcal{N})$, so it suffices to
show that $S$ is of continuum cardinality. But this follows from a
well-known fact of classical DST (Theorem 29.1 in \cite{ke95}).
 \qed

\begin{remark}\label{contin1}
The last result  can not be improved within ZFC because, as is well
known, it is consistent with  ZFC that there is a non-countable set
$S\in\mathbf{\Pi}^1_1(\mathcal{N})$ of cardinality less than
continuum.
\end{remark}

Next we establish an extension of the Suslin theorem which equates
the Borel sets to the $\mathbf{\Delta}^1_1$-sets. This is a
classical result of DST for the case of Polish spaces, and it was
extended to quasi-Polish spaces in \cite{br}.

For this we need the following version of a well-known easy fact:

\begin{lemma}\label{subsp}
Let $X\subseteq Y$ be topological spaces and $1\leq\alpha<\omega_1$.
Then $\mathbf{\Sigma}^1_\alpha(X)=\{X\cap A\mid
A\in\mathbf{\Sigma}^1_\alpha(Y)\}$,
$\mathbf{\Pi}^1_\alpha(X)=\{X\cap A\mid
A\in\mathbf{\Pi}^1_\alpha(Y)\}$, and similarly for the Borel
hierarchy.
\end{lemma}

\begin{proposition}\label{suslin}
The Suslin theorem holds for any space $X$ in
$\CBZ(\mathbf{\Delta}^1_1)$, i.e.
$\mathbf{\Delta}^1_1(X)=\bigcup\{\mathbf{\Sigma}^0_\alpha(X)\mid
\alpha<\omega_1\}$.
\end{proposition}

{\bf Proof.} It suffices to show the inclusion
$\mathbf{\Delta}^1_1(X)\subseteq\bigcup\{\mathbf{\Sigma}^0_\alpha(X)\mid
\alpha<\omega_1\}$. Assume without loss of generality that
$X\in\mathbf{\Delta}^1_1(P\omega)$. Let
$A\in\mathbf{\Delta}^1_1(X)$, then $A$ is in both
$\mathbf{\Sigma}^1_1(X)$ and $\mathbf{\Pi}^1_1(X)$. By Lemma
\ref{subsp}, $A=X\cap B=X\cap C$ for some
$B\in\mathbf{\Sigma}^1_1(P\omega)$ and
$C\in\mathbf{\Pi}^1_1(P\omega)$. Then
$A\in\mathbf{\Delta}^1_1(P\omega)$. By Suslin theorem for $P\omega$,
$A\in\mathbf{\Sigma}^0_\alpha(P\omega)$ for some $\alpha<\omega_1$.
By Lemma \ref{subsp}, $A\in\mathbf{\Sigma}^0_\alpha(X)$.
 \qed

As is well known \cite{ke95}, the Borel and Luzin hierarchies do not
collapse in any Polish uncountable space $X$ (for the Borel
hierarchy, for instance, this means that
$\mathbf{\Sigma}^0_\alpha(X)\not=\mathbf{\Pi}^0_\alpha(X)$ for any
$\alpha<\omega_1$). In \cite{br} this was extended to the
quasi-Polish spaces (which coincide with the spaces in
$\CBZ(\mathbf\Pi^0_2)$). We conclude this section with a further
extension of the non-collapse property. For this we need the
following lemma:

\begin{lemma}\label{noncol}
Let $1\leq\alpha<\omega_1$, $X ,Y$ be
$\mathbf{\Sigma}^0_\alpha$-isomorphic topological spaces, and the
Borel hierarchy (resp. the hyperprojective hierarchy) in $X$ does
not collapse. Then the Borel hierarchy (resp. the hyperprojective
hierarchy) in $Y$ does not collapse.
\end{lemma}

{\bf Proof.}  The both hierarchies are treated similarly, so
consider only the Borel hierarchy. Suppose for a contradiction that
$\mathbf{\Sigma}^0_\beta(Y)=\mathbf{\Pi}^0_\beta(Y)$ for some
$\beta<\omega_1$. By the definition of the Borel hierarchy,
$\mathbf{\Sigma}^0_\gamma(Y)=\mathbf{\Pi}^0_\gamma(Y)$ for all
countable ordinals $\gamma\geq\beta$, in particular for
$\gamma=sup\{\alpha,\beta\}$. By Lemma \ref{gfam1},
$\mathbf{\Sigma}^0_\gamma(X)=\mathbf{\Pi}^0_\gamma(X)$. A
contradiction.
 \qed

\begin{proposition}\label{noncol1}
The Borel  and hyperprojective hierarchies do not collapse for any
uncountable space $X$ in $\CBZ(\mathbf{\Delta}^1_1)$.
\end{proposition}

{\bf Proof.}  By Proposition \ref{aliso} and Lemma \ref{noncol}, we
can without loss of generality assume that $X$ is a subspace of
$\mathcal{N}$ such that $X\in\mathbf{\Delta}^1_1(\mathcal{N})$. By
Theorem 29.1 in \cite{ke95}, there is a subspace $C\subseteq X$
homeomorphic to the Cantor space. Since the Borel and
hyperprojective hierarchies in $C$ do not collapse, by Lemma
\ref{subsp} they also do not collapse in $X.$
 \qed

%%%%%%%%%%%%%%%%%%%%%%%%%%%%%%%%%%%%%%%%%%%%%%%%%%%%%%%

\section{Difference hierarchies}\label{opens}

In this section we extend some classical facts on the Hausdorff
difference hierarchy (DH) of sets, like the Hausdorff-Kuratowski
theorem and the non-collapse property in Polish spaces, to larger
classes of cb$_0$-spaces and to the case of $k$-partitions.

\subsection{Difference hierarchies of sets}\label{dhsets}

Here  we recall definition and  basic properties of the Hausdorff
difference hierarchy of sets, and extend some facts on the DH in
Polish spaces to larger classes of cb$_0$-spaces.

An ordinal $\alpha$ is {\em even}  (resp.
{\em odd}) if $\alpha=\lambda+n$ where $\lambda$ is either zero or
a limit ordinal and $n<\omega$, and the number $n$ is even (resp.,
odd). For an ordinal $\alpha$, let $r(\alpha)=0$ if $\alpha$ is
even and $r(\alpha)=1$, otherwise. For any ordinal $\alpha$,
define the operation $D_\alpha$ sending sequences of sets
$\{A_\beta\}_{\beta<\alpha}$ to sets by
 $$D_\alpha(\{A_\beta\}_{\beta<\alpha})=
\bigcup\{A_\beta\setminus\bigcup_{\gamma<\beta}A_\gamma\mid
\beta<\alpha,\,r(\beta)\not=r(\alpha)\}.$$

For any ordinal $\alpha<\omega_1$ and  any pointclass ${\cal L}$ in
$X$, let $D_\alpha({\cal L})$ be the class of all sets
$D_\alpha(\{A_\beta\}_{\beta<\alpha})$, where $A_\beta\in{\cal L}$
for all $\beta<\alpha$. By the {\em difference hierarchy over ${\cal
L}$} we mean the sequence $\{D_\alpha({\cal
L})\}_{\alpha<\omega_1}$. Usually we assume that ${\cal L}$ is a
{\em base} which by definition means that ${\cal L}$ is closed under
finite intersection and countable union (note that in finitary
versions of the DH we used the term ``$\sigma$-base'' to denote such
pointclasses but, since we are interested here only in such
pointclasses, we simplify the terminology). As usual, classes
$D_\alpha({\cal L})$, $\check{D}_\alpha({\cal L})$ are called
non-self-dual levels while $D_\alpha({\cal
L})\cap\check{D}_\alpha({\cal L})$ are called self-dual levels of
the DH.

Over bases, the difference hierarchy really looks as a hierarchy,
i.e., any level and its dual are contained in all higher levels. The
most interesting cases for Descriptive Set Theory are difference
hierarchies over non-zero levels of the Borel hierarchy, whose
$\bf{\Sigma}$-levels are
${\bf\Sigma}^{-1,\theta}_\alpha(X)=D_\alpha({\bf\Sigma}^0_\theta(X))$,
for any space $X$  and for all $\alpha,\theta<\omega$, $\theta>0$.
For $\theta=1$, we simplify ${\bf\Sigma}^{-1,\theta}_\alpha$ to
${\bf\Sigma}^{-1}_\alpha$.

A classical result of DST is the following Hausdorff-Kuratowski theorem:

\begin{theorem}\label{hk}
Let $X$ be a Polish space. For any non-zero ordinal $\theta<\omega_1$,
 $\bigcup\{{\bf\Sigma}^{-1,\theta}_\alpha(X)\mid \alpha<\omega_1\}=
 \mathbf\Delta^0_{\theta+1}(X).$
\end{theorem}

In \cite{br} this result was extended to the quasi-Polish spaces.
This extension is an easy corollary of the following nice result
(Theorem 68 in \cite{br} based on the important Lemma 17 in \cite{sr07}):

\begin{theorem}\label{hk1}
Let  $X$ be a cb$_0$-space, $\delta:D\to X$ an admissible
representation of $X$ ($D\subseteq\mathcal{N}$), $A\subseteq X$ $X$,
$\alpha,\theta<\omega_1$  and $\theta\geq1$. Then $A\in
D_\alpha(\mathbf{\Sigma}^0_\theta(X))$ iff $\delta^{-1}(A)\in
D_\alpha(\mathbf{\Sigma}^0_\theta(D))$.
\end{theorem}

One  of the aims of this section is to extend these result from
sets to $k$-partitions. This needs some information on $k$-forests
and $h$-preorders which are recalled in the next subsection.

Next  we establish a partial extension of the Hausdorff-Kuratowski
theorem for Polish spaces (which coincide with the spaces in
$\CBZ(\mathbf\Pi^0_2)$) to a larger class of cb$_0$-spaces:

\begin{proposition}\label{hkext}
For  any space $X$ in $\CBZ(\mathbf{\Delta}^1_1)$ there is a
non-zero ordinal $\beta<\omega_1$ such that the Hausdorff-Kuratowski
theorem holds in $X$ for each countable ordinal $\theta\geq\beta$,
i.e. $\bigcup\{{\bf\Sigma}^{-1,\theta}_\alpha(X)\mid
\alpha<\omega_1\}=\mathbf\Delta^0_{\theta+1}(X)$.
\end{proposition}

{\bf Proof.}  Assume without loss of generality that
$X\in\mathbf{\Delta}^1_1(P\omega)$. By Suslin theorem for $P\omega$,
$X\in\mathbf{\Sigma}^0_\beta(P\omega)$ for some $\beta<\omega_1$. It
remains to show that
$\mathbf\Delta^0_{\theta+1}(X)\subseteq\bigcup\{{\bf\Sigma}^{-1,\theta}_\alpha(X)\mid
\alpha<\omega_1\}$ for all $\theta\geq\beta$. Let
$A\in\mathbf\Delta^0_{\theta+1}(X)$, then $A$ is in both
$\mathbf\Sigma^0_{\theta+1}(X)$ and $\mathbf\Pi^0_{\theta+1}(X)$. By
Lemma \ref{subsp}, $A=X\cap B=X\cap C$ for some
$B\in\mathbf\Sigma^0_{\theta+1}(P\omega)$ and
$C\in\mathbf\Pi^0_{\theta+1}(P\omega)$. Then
$A\in\mathbf\Delta^0_{\theta+1}(P\omega)$.  By the
Hausdorff-Kuratowski theorem for $P\omega$,
$A\in\mathbf{\Sigma}^{-1,\theta}_\alpha(P\omega)$ for some
$\alpha<\omega_1$. Since $\theta\geq\beta$ and
$X\in\mathbf{\Sigma}^0_\beta(P\omega)$,
$A\in\mathbf{\Sigma}^{-1,\theta}_\alpha(X)$.
 \qed

The problem of the non-collapse of the DHs is more subtle (compared
with the problem of non-collapse of the Borel and Luzin hierarchies)
but it is again possible to prove at least a partial result about
this property. First we formulate an analogue of Lemma \ref{noncol}
which is proved essentially by the same argument:

\begin{lemma}\label{noncoldh}
Let $1\leq\alpha<\omega_1$, $X ,Y$ be
$\mathbf{\Sigma}^0_\alpha$-isomorphic topological spaces,
$\alpha\leq\theta<\omega_1$, and the DH
$\{\mathbf{\Sigma}^{-1,\theta}_\beta(X)\}_{\beta<\omega_1}$ does not
collapse. Then the DH
$\{\mathbf{\Sigma}^{-1,\theta}_\beta(Y)\}_{\beta<\omega_1}$ does not
collapse.
\end{lemma}

Once we have this  lemma and note that Lemma \ref{subsp} holds also
for the DHs, we easily deduce the following:

\begin{proposition}\label{noncoldh1}
For any  uncountable space $X$ in $\CBZ(\mathbf{\Delta}^1_1)$ there
is a non-zero ordinal $\alpha<\omega_1$ such that the DH
$\{\mathbf{\Sigma}^{-1,\theta}_\beta(X)\}_{\beta<\omega_1}$ does not
collapse for each countable ordinal $\theta\geq\alpha$.
\end{proposition}

\subsection{h-Preorder}\label{hpre}

Here we discuss some posets which serve as notation systems for
levels of the DHs of $k$-partitions. All notions and facts of this
subsection are contained (at least, implicitly) in \cite{s07,s07a}.

Posets considered here are assumed to be (at most) countable and
without infinite chains. The absence of infinite chains in a poset
$(P;\leq)$ is of course equivalent to well-foundednes of both
$(P;\leq)$ and $(P;\geq)$. By a {\em  forest} we mean a poset
without infinite chains in which every upper cone $\{y\mid x\leq
y\}$ is a chain. {\em  Tree} is a forest having the biggest element
(called {\em the root} of the tree).

Let $(T;\leq)$ be a tree without infinite chains; in particular, it
is well-founded. As for each well-founded partial order, there is a
canonical rank function $rk_T$ from $T$ to ordinals.  The rank
$rk(T)$ of  $(T;\leq)$ is by definition the ordinal $rk_T(r)$, where
$r$ is the root of $(T;\leq)$.  It is well-known that the rank of
any countable tree without infinite chains is a countable ordinal,
and any countable ordinal is the rank of such a tree.

A {\em $k$-poset}  is a triple $(P;\leq,c)$ consisting of a poset
$(P;\leq)$ and a labeling $c:P\rightarrow k$. Rank of a $k$-poset
$(T;\leq,c)$ is by definition the rank of $(T;\leq)$. {\em Morphism}
$f:(P;\leq,c)\rightarrow(P^\prime;\leq^\prime,c^\prime)$ between
$k$-posets is a monotone function
$f:(P;\leq)\rightarrow(P^\prime;\leq^\prime)$ respecting the
labelings, i.e. satisfying $c=c^\prime\circ f$. Let
 $\widetilde{\mathcal F}_k$ and
$\widetilde{\mathcal T}_k$  denote the classes of all  countable
$k$-forests and countable $k$-trees  without infinite chains,
respectively. Note that we use tilde in our notation in order to
distinguish the introduced objects from their finitary versions
extensively studied in \cite{s06,s08a,s12}.

The $h$-preorder $\leq_h$ on $\widetilde{\mathcal P}_k$ is defined
as follows: $(P,\leq,c)\leq_h(P^\prime,\leq^\prime,c^\prime)$, if
there is a morphism from $(P,\leq,c)$ to
$(P^\prime,\leq^\prime,c^\prime)$.  Let  $\widetilde{\mathbb
F}_k,\widetilde{\mathbb T}_k$ be the quotient posets of
$\widetilde{\mathcal F}_k$ and $\widetilde{\mathcal T}_k$ under the
$h$-equivalence, respectively. Let $\widetilde{\mathbb F}_k^\prime$
be obtained from $\widetilde{\mathbb F}_k$ by adjoining a new
smallest element $\bot$ (corresponding to the empty forest).

Let $P\sqcup Q$ be the disjoint union of $k$-posets $P,Q$ and $\bigsqcup_i
P_i=P_0\sqcup P_1\sqcup\cdots$  the disjoint union of an infinite sequence
$P_0,P_1,\ldots$ of $k$-posets. For a $k$-forest $F$ and  $i<k$,
 let $p_i(F)$ be the $k$-tree
obtained from $F$ by adjoining a new biggest element and assigning
the label $i$ to this element. It is clear that the introduced
operations respect the  $h$-equivalence and that any countable
$k$-forest is $h$-equivalent to a countable term of signature
$\{\sqcup,p_0,\ldots,p_{k-1},0,\ldots,k-1\}$ without free variables
(the constant symbol $i$ in the signature is interpreted as the
singleton tree carrying the label $i$).

\begin{proposition}\label{main2}
\begin{enumerate}
 \item For any $k\geq 2$, the structures $(\widetilde{\mathcal F}_k;\leq)$
and $(\widetilde{\mathcal T}_k;\leq)$ are well preorders of rank
$\omega_1$.
 \item The posets $\widetilde{\mathbb F}_2$ and
 $\widetilde{\mathbb T}_2$ have width 2 (i.e., they have
no antichains with more that 2 elements).
 \item The poset $\widetilde{\mathbb F}_k^\prime$ is a distributive lattice
where any countable set of elements have a supremum.
 \item The set $\sigma ji(\widetilde{\mathbb F}_k^\prime)$  of
 $\sigma$-join-irreducible elements   of the
 lattice $\widetilde{\mathbb F}_k^\prime$ (i.e., the elements $x$ such that
 $x\leq\bigsqcup\{y_n\mid n<\omega\}$ implies that $y\leq y_n$ for some
 $n$) coincides with $\widetilde{\mathbb T}_k$.
 \item The set $ ji(\widetilde{\mathbb F}_k^\prime)$ of join-irreducible elements of the
 lattice $\widetilde{\mathbb F}_k^\prime$ coincides with $\widetilde{\mathbb T}_k\cup S$
 where $S$ is the set of supremums of infinite increasing sequences of
 elements in $\widetilde{\mathbb T}_k$.
 \item Any element of $\widetilde{\mathbb F}_k^\prime$ is the infimum of finitely many
 elements of $\widetilde{\mathbb T}_k$.
\end{enumerate}
\end{proposition}

For a result  in the next subsection we need the following canonical
representatives for the structures $\widetilde{\mathcal F}_2$ and
$\widetilde{\mathcal T}_2$. Define by induction the sequence
$\{T_\alpha\}_{\alpha<\omega_1}$ as follows: $T_0=0$,
$T_{\alpha+1}=p_0(\overline{T}_\alpha)$ where $\overline{T}_\alpha$
is obtained from $T_\alpha$ by changing any label $l<2$ by the label
$1-l$, and
$T_\lambda=p_0(\overline{T}_{\alpha_0}\sqcup\overline{T}_{\alpha_1}\sqcup\cdots)$
for a limit ordinal $\lambda$ where $\alpha_0<\alpha_1<\cdots$is a
sequence of odd ordinals satisfying $sup\{\alpha_n\mid
n<\omega\}=\lambda$. The next assertion   follows from the proof of
the corresponding result in \cite{s07,s07a}.

\begin{proposition}\label{t2}
\begin{enumerate}
 \item For all $\alpha<\beta<\omega_1$, $T_\alpha\sqcup\overline{T}_\alpha<_hT_\beta$
and $T_\alpha,\overline{T}_\alpha$ are $\leq_h$-incomparable.
 \item  Any element of $\widetilde{\mathcal T}_2$ (resp. of $\widetilde{\mathcal F}_2$)
is $h$-equivalent to precisely one of
$T_\alpha,\overline{T}_\alpha$, (resp. to precisely one of
$T_\alpha,\overline{T}_\alpha,
T_\alpha\sqcup\overline{T}_\alpha,\bigsqcup_{\alpha<\lambda}T_\alpha$
where $\lambda$ is a limit countable ordinal).
 \item   $T_\alpha\sqcup\overline{T}_\alpha\equiv_hT_{\alpha+1}
\sqcap\overline{T}_{\alpha+1}$ and
$\bigsqcup_{\alpha<\lambda}T_\alpha\equiv_hT_{\lambda}\sqcap\overline{T}_{\lambda}$
for each limit countable ordinal $\lambda$, where $\sqcap$ is the
binary operation of $k$-forests inducing the infimum operation in
$\widetilde{\mathbb F}_k^\prime$.
\end{enumerate}
\end{proposition}

\subsection{Difference hierarchies of $k$-partitions}\label{kpart}

Here we extend the difference hierarchy of sets to that of
$k$-partitions. Note that similar hierarchies were considered e.g.
in \cite{s07,s07a} but in fact the definition here is slightly
different which results in equivalence of the corresponding DHs only
over the bases with $\omega$-reduction properties. In general, the
definition here is better than  the previous one because it yields,
for example, Theorem \ref{dhpart6} for $\theta=1$ which is not
always the case for the previous definition.

Let ${\cal L}$ be a  base in $X$. Recall that a $k$-partition of $X$
is a function $A:X\rightarrow k$ often written as a tuple
$(A_0,\ldots,A_{k-1})$ where $A_i=\{x\in X\mid A(x)=i\}$. By a {\em
partial $k$-partition of $X$} we mean a function $A:Y\rightarrow k$
for some $Y\in\mathcal{L}$. Let $P\in\widetilde{\mathcal{F}}_k$. We
say that a partial $k$-partition $A$ is defined by a $P$-family
$\{B_p\}_{p\in P}$ of $\mathcal{L}$-sets if $A_i=\bigcup_{p\in
P_i}\tilde{B}_p$ for each $i<k$ where
$\tilde{B}_p=B_p\setminus\bigcup_{q<p}B_q$ and $P_i=c^{-1}(i)$; note
that in this case $A\in k^Y$ where $Y=\bigcup_{p\in P}B_p$.

We denote by ${\mathcal L}^Y(P)$ the set of partitions
$A:Y\rightarrow k$ defined by some $P$-family $\{B_p\}_{p\in P}$ of
$\mathcal{L}$-sets. In case $Y=X$ we omit the superscript $X$ and
call the classes ${\mathcal L}(P),P\in\widetilde{\mathcal{F}}_k$
{\em  levels of the DH of $k$-partitions over $\mathcal{L}$}. We
formulate some basic properties of the levels. We omit the proofs
because they are quite similar to the corresponding proofs for the
finitary version of the DH in \cite{s12}.

\begin{proposition}\label{treedh1}
 \begin{enumerate}\itemsep-1mm
 \item  If $A\in{\mathcal L}^Y(P)$ then $A|_Z\in{\mathcal L}^Z(P)$
for each $Z\subseteq Y$, $Z\in\mathcal{L}$.
 \item  Any $A\in{\mathcal L}^Y(P)$ is defined by a monotone
 $P$-family $\{C_p\}$ (monotonicity means that $C_q\subseteq C_p$ for
 $q\leq p$).
 \item Let $f$ be a function on $X$ such that $f^{-1}(A)\in\mathcal{L}$
for each $A\in\mathcal{L}$. Then $A\in{\mathcal L}^Y(P)$ implies
$f^{-1}(A)=(f^{-1}(A_0),\ldots,f^{-1}(A_{k-1}))\in{\mathcal
L}^{f^{-1}(Y)}(P)$.
 \item If $P\leq_hQ$ then ${\mathcal
L}^Y(P)\subseteq{\mathcal L}^Y(Q)$.
 \item For all $F,G\in\widetilde{\mathcal{F}}_k$, $\mathcal{L}(F)\cap
 \mathcal{L}(G)=\mathcal{L}(F\sqcap G)$.
 \end{enumerate}
  \end{proposition}

Proposition \ref{treedh1}(5) and Proposition \ref{main2}(6) show
that any level $\mathcal{L}(F)$ is a finite intersection of the
levels $\mathcal{L}(T)$, $T\in\widetilde{\mathcal{T}}_k$. This
remark, together with the results below, suggest that the levels
$\mathcal{L}(T)$, $T\in\widetilde{\mathcal{T}}_k$ are analogs for
the DH of $k$-partitions of the non-self-dual levels
$D_\alpha(\mathcal{L}),\check{D}_\alpha(\mathcal{L})$ of the DH of
sets. Therefore, the precise analog of the DH of sets is the family
$\{\mathcal{L}(T)\}_{T\in\widetilde{\mathcal{T}}_k}$ rather than the
family $\{\mathcal{L}(F)\}_{F\in\widetilde{\mathcal{F}}_k}$.

The meaning  of the last paragraph might be not clear because it is
not even obvious that the DH of $k$-partitions really extends the DH
of sets; we have at least to show that the DH of $2$-partitions
essentially coincides with the DH of sets. We do this in the next
proposition where we employ the 2-trees $T_\alpha$ from Subsection
\ref{hpre}.

\begin{proposition}\label{main3}
Let  $\mathcal{L}$ be a base in $X$. Then
$\mathcal{L}(T_\alpha)=D_\alpha(\mathcal{L})$ for each
$\alpha<\omega_1$.
\end{proposition}

{\bf Proof.} Let $A\in\mathcal{L}(T_\alpha)$ be defined by a family
$\{B_p\}_{p\in T_\alpha}$ where $B_p\in\mathcal{L}$. Define the
sequence $\{A_\beta\}_{\beta<\alpha}$ as follows: if
$r(\beta)=r(\alpha)$ then $A_\beta=\bigcup\{B_p\mid
rk(p)\leq\beta\wedge c(p)=0\}$, otherwise $A_\beta=\bigcup\{B_p\mid
rk(p)\leq\beta\wedge c(p)=1\}$. Then $A_\beta\in\mathcal{L}$ and
$A=\{\tilde{A}_\beta\mid r(\beta)\not=r(\alpha)\}$ where
$\tilde{A}_\beta=A_\beta\setminus\bigcup_{\gamma<\beta}A_\gamma$,
hence $A\in D_\alpha(\mathcal{L})$.

Conversely,  let $A\in D_\alpha(\mathcal{L})$, then
$A=\{\tilde{A}_\beta\mid r(\beta)\not=r(\alpha)\}$ for some sequence
$\{A_\beta\}_{\beta<\alpha}$ of $\mathcal{L}$-sets. Define the
family $\{B_p\}_{p\in T_\alpha}$ of $\mathcal{L}$-sets as follows:
if $p$ is the root of $T_\alpha$ then $B_p=X$, otherwise
$B_p=A_{rk(p)}$ where $rk:T_\alpha\to\alpha+1$ is the rank function.
Since $rk$ is a surjection for each $\alpha$,
$\tilde{B}_p=\tilde{A}_{rk(p)}$.

Note that  $\tilde{B}_p,\tilde{B}_q$ are disjoint whenever
$c(p)\not=c(q)$ because $\tilde{A}_\beta,\tilde{A}_\gamma$ are
disjoint for distinct $\beta,\gamma<\alpha$. If $p$ is the root of
$T_\alpha$ then
$\tilde{B}_p=X\setminus\bigcup_{\beta<\alpha}\tilde{A}_\beta\subseteq\overline{A}$.
If $p$ is not the root then, by the definition of $T_\alpha$,
$r(rk(p))=r(\alpha)$ iff $c(p)=0$. Then for each
$x\in\bigcup_{\beta<\alpha}\tilde{A}_\beta$ we have: $x\in A$ iff
$r(\beta)\not=r(\alpha)$ (where $\beta$ is the unique ordinal with
$x\in\tilde{A}_\beta$) iff $c(p)=1$. Thus,
$A\in\mathcal{L}(T_\alpha)$ is defined by  $\{B_p\}_{p\in T_\alpha}$
and therefore $A\in\mathcal{L}(T_\alpha)$.
 \qed

For $F\in\widetilde{\mathcal{F}}_k$, by a {\em reduced $F$-family of
$\mathcal{L}$-sets} we mean a monotone $F$-family $\{B_p\}_{p\in F}$
of $\mathcal{L}$-sets such that $B_p\cap B_q=\emptyset$ for all
incomparable $p,q\in F$. Let $A\in\mathcal{L}^Y_r(F)$ be the set of
partial $k$-partitions defined by reduced $F$-families
$\{B_p\}_{p\in F}$ of $\mathcal{L}$-sets such that $\bigcup_pB_p=Y$.
The next result is an infinitary version of Proposition 7.15
\cite{s12} and is proved by essentially the same argument.

\begin{proposition}\label{dhpart2}
Let $\mathcal{L}$ have the $\omega$-reduction property, $Y\in\mathcal{L}$
and $F\in\mathcal{F}_k$. Then
$\mathcal{L}^Y(F)=\mathcal{L}^Y_r(F)$, in particular $\mathcal{L}(F)=\mathcal{L}_r(F)$.
\end{proposition}

In \cite{s13},  principal total representations (TR) of pointclasses
were introduced and studied. The results in \cite{s13} naturally
extend to $k$-partitions. By a {\em family of $k$-partition classes}
we mean a family $ \mathbf{\Gamma}=\{\mathbf{\Gamma}(X)\} $ indexed
by arbitrary topological spaces $X$ such that  $
\mathbf{\Gamma}(X)\subseteq k^X$ and $A\circ f\in\mathbf{\Gamma}(X)$
for any continuous function $f:X\to Y$ and any $k$-partition $A:Y\to
k$ from $\mathbf{\Gamma}(Y)$.

We note that, by Proposition \ref{treedh1}(3),  the levels of the
DHs over any $\mathbf{\Sigma}$-level of the Borel or Luzin hierarchy
are families of $k$-partition classes. Let
$\mathbf{\Sigma}^0_\theta(X,F)$ be the $F$-level
($F\in\widetilde{\mathcal{F}}_k$) of the DH of $k$-partitions over
$\mathbf{\Sigma}^0_\theta(X)$. The next fact is obvious:

\begin{proposition}\label{dhpart3}
Let  $\theta$ be a non-zero countable ordinal and
$F\in\widetilde{\mathcal{F}}_k$. Then
$\mathbf{\Sigma}^0_\theta(F)=\{\mathbf{\Sigma}^0_\theta(X,F)\}_X$ is
a family of $k$-partition classes.
\end{proposition}

Let  $\{\mathbf{\Gamma}(X)\} $ be a family of $k$-partition classes.
A function $\nu:\mathcal{N}\to\mathbf{\Gamma}(X)$ is a {\em
$\mathbf{\Gamma}$-TR} if the $k$-partition $\lambda a,x.\nu(a)(x)$
is in $\mathbf{\Gamma}(\mathcal{N}\times X)$. Such $\nu$ is a {\em
principal $\mathbf{\Gamma}$-TR} if any $\mathbf{\Gamma}$-TR
$\mu:\mathcal{N}\to\mathbf{\Gamma}(X)$ is continuously reducible to
$\nu$.

According  to Theorem 5.2 in \cite{s13}, the non-self-dual levels
$\mathbf{\Gamma}$-TR of the classical hierarchies in cb$_0$-spaces
have $\mathbf{\Gamma}(X)$-TRs. This result extends to the
``non-self-dual'' levels of the DHs of $k$-partitions but only under
the additional assumption that the corresponding bases have the
$\omega$-reduction property:

\begin{proposition}\label{dhpart4}
Let  $X$ be a cb$_0$-space, $\theta\geq2$  a  countable ordinal and
$T\in\widetilde{\mathcal{T}}_k$. Then
$\mathbf{\Sigma}^0_\theta(X,T)$ has a principal
$\mathbf{\Sigma}^0_\theta(T)$-TR. If $X$ is in addition
zero-dimensional then $\mathbf{\Sigma}^0_1(X,T)$ has a principal
$\mathbf{\Sigma}^0_1(T)$-TR.
\end{proposition}

{\bf Proof Hint.}  Modulo Theorem 5.2 in \cite{s13}, the proof is
straightforward, so we give only informal proof hints. Note that, by
Proposition \ref{uniform}, the assumptions on $\theta$ guarantee
that the class $\mathbf{\Sigma}^0_\theta(X)$ has the
$\omega$-reduction property, so by Proposition \ref{dhpart2} the
elements of $\mathbf{\Sigma}^0_\theta(X,T)$ are precisely those
defined by the monotone reduced $T$-families $\{B_p\}$ of
$\mathbf{\Sigma}^0_\theta(X)$-sets with $B_p=X$ where $p$ is the
root of $T$. The principal $\mathbf{\Sigma}^0_\theta$-TR of
$\mathbf{\Sigma}^0_\theta(X)$ induces a natural representaion of all
$T$-families $\{C_p\}$ of $\mathbf{\Sigma}^0_\theta(X)$-sets with
$C_p=X$. The problem is that such a family typically does not define
any $k$-partition, hence we do not in general have an induced TR of
$\mathbf{\Sigma}^0_\theta(X,T)$. But the $\omega$-reduction property
gives a uniform procedure of transforming $\{C_p\}$ to a monotone
reduced $T$-family $\{B_p\}$ of $\mathbf{\Sigma}^0_\theta(X)$-sets
with $B_p=X$, such that $\{B_p\}=\{C_p\}$ whenever $\{C_p\}$ already
has this property. Since any such $\{B_p\}$ defines an element of
$\mathbf{\Sigma}^0_\theta(X,T)$, this  induces a TR of
$\mathbf{\Sigma}^0_\theta(X,T)$. It is straightforward to check that
this TR has the desired properties. \qed

We conclude  this section with extending the Hausdorff-Kuratowski
theorem to $k$-partitions. First we extend Theorem \ref{hk1} to
$k$-partitions:

\begin{theorem}\label{dhpart5}
Let $X$ be  a cb$_0$-space, $\delta:D\to X$ an admissible
representation of $X$ ($D\subseteq\mathcal{N}$), $A:X\to k$ a
$k$-partition of $X$, $\theta\geq1$  a  countable ordinal and
$F\in\widetilde{\mathcal{F}}_k$. Then
$A\in\mathbf{\Sigma}^0_\theta(X,F)$ iff
$A\circ\delta\in\mathbf{\Sigma}^0_\theta(D,F)$.
\end{theorem}

{\bf Proof}  is similar to the proof of Theorem 68 in \cite{br}. Let
first $A\in\mathbf{\Sigma}^0_\theta(X,T)$, then $A$ is defined by an
$F$-family $\{B_p\}$ of $\mathbf{\Sigma}^0_\theta(X)$-sets. Then
$A\circ\delta$ is defined by the $F$-family $\{\delta^{-1}(B_p)\}$
of $\mathbf{\Sigma}^0_\theta(D)$-sets, hence
$A\circ\delta\in\mathbf{\Sigma}^0_\theta(D,T)$.

Conversely,  let $A\circ\delta\in\mathbf{\Sigma}^0_\theta(D,T)$,
then $A\circ\delta$ is defined by an $F$-family $\{C_p\}$ of
$\mathbf{\Sigma}^0_\theta(D)$-sets, so
$\delta^{-1}(A_i)=\bigcup\{\tilde{C}_p\mid c(p)=i\}$ for each $i<k$.
By the proof of Theorem 68 in \cite{br}, we can without loss of
generality assume that $\delta$ has Polish fibers, i.e.,
$\delta^{-1}(x)$ is Polish for each $x\in X$. For any $p\in F$, let
$B_p$ consist of the elements $x\in X$ such that the set
$C_p\cap\delta^{-1}(x)$ is non-meager in $\delta^{-1}(x)$. By the
proof of Theorem 68 in \cite{br},
$B_p\in\mathbf{\Sigma}^0_\theta(X)$, hence it suffices to show that
$A$ is defined by the $F$-family $\{B_p\}$.

First we check that $\tilde{B}_p\subseteq\delta(\tilde{C}_p)$. Let
$x\in\tilde{B}_p$, so $C_p\cap\delta^{-1}(x)$ is non-meager in
$\delta^{-1}(x)$ and $C_q\cap\delta^{-1}(x)$ is meager in
$\delta^{-1}(x)$ for each $q<p$, hence also
$(\bigcup_{q<p}C_q)\cap\delta^{-1}(x)$ is meager in
$\delta^{-1}(x)$. Since
 $$C_p\cap\delta^{-1}(x)=(\tilde{C}_p\cap\delta^{-1}(x))\cup
(\bigcup_{q<p}C_q)\cap\delta^{-1}(x),$$
$\tilde{C}_p\cap\delta^{-1}(x)$ is non-meager in $\delta^{-1}(x)$,
in particular $\tilde{C}_p\cap\delta^{-1}(x)$ is non-empty. Let
$a\in\tilde{C}_p\cap\delta^{-1}(x)$, then
$x=\delta(a)\in\delta(\tilde{C}_p)$, as desired.

We  have to show that $A_i=\bigcup\{\tilde{B}_p\mid c(p)=i\}$ for
each $i<k$. Let first $x\in\tilde{B}_p$, $c(p)=i$. Then
$x=\delta(a)$ for some $a\in\tilde{C}_p\subseteq \delta^{-1}(A_i)$.
Thus, $x\in A_i$.

Conversely, let $x\in A_i$. Choose $a\in D$ with $x=\delta(a)$. Then
$a\in\delta^{-1}(A_i)$, so $a\in\tilde{C}_p$ for some $p\in F$,
$c(p)=i$. Note that $x\in B_q$ for some $q\in F$ (otherwise,
$C_q\cap\delta^{-1}(x)$ is meager in $\delta^{-1}(x)$ for each $q\in
F$, hence also $\delta^{-1}(x)$ is meager, a contradiction). Then
$x\in\tilde{B}_q$ for some $q\in F$, hence $x=\delta(b)$ for some
$b\in\tilde{C}_q$. Let $j=c(q)$, then $a\in\delta^{-1}(A_i)$ and
$b\in\delta^{-1}(A_j)$, then $x\in A_i\cap A_j$, so $i=j$ and
$c(q)=i$. \qed

As an immediate  corollary, we obtain the Hausdorff-Kuratowski
theorem for $k$-partitions in quasi-Polish spaces.

\begin{theorem}\label{dhpart6}
Let $X$  be a quasi-Polish space and $\theta\geq1$  a  countable
ordinal. Then $\bigcup\{\mathbf{\Sigma}^0_\theta(X,T)\mid
F\in\widetilde{\mathcal{F}}_k\}=(\mathbf{\Delta}^0_{\theta+1}(X))_k$.
\end{theorem}

{\bf Proof.} For $\theta\geq2$, the assertion follows from Theorem
5.1 in \cite{s08} and Theorem \ref{hk1} but for $\theta=1$ the
result is new. Let $A\in(\mathbf{\Delta}^0_2(X))_k$, we have to show
that $A\in\bigcup\{\mathbf{\Sigma}^0_1(X,T)\mid
F\in\widetilde{\mathcal{F}}_k\}$. Let $\delta$ be an admissible
total representation of $X$. Then $A\circ\delta
\in(\mathbf{\Delta}^0_2(\mathcal{N}))_k$. Since
$\mathbf{\Sigma}^0_1(\mathcal{N})$ has the $\omega$-reduction
property,
$A\circ\delta\in\bigcup\{\mathbf{\Sigma}^0_\alpha(\mathcal{N},F)\mid
F\in\widetilde{\mathcal{F}}_k\}$ by Theorem 5.1 in \cite{s08}. By
Theorem \ref{dhpart5}, $A\in\bigcup\{\mathbf{\Sigma}^0_1(X,T)\mid
F\in\widetilde{\mathcal{F}}_k\}$.
 \qed

\begin{remark}
Using the argument of Proposition \ref{hkext}, it is straightforward
to extend the last result to spaces $X$ in
$\CBZ(\mathbf{\Delta}^1_1)$.
 \end{remark}

%%%%%%%%%%%%%%%%%%%%%%%%%%%%%%%%%%%%%%%%%%%%%%%%%%%%%%%%%%%%
%
%
%%%%%%%%%%%%%%%%%%%%%%%%%%%%%%%%%%%%%%%%%%%%%%%%%%%%%%%%%%%%

\section{Fine hierarchies}\label{fhs}

In this section we extend the DH of $k$-partitions to the FH of
$k$-partitions.  Many results and proofs here extend the  ones for
the DH from the previous section or the corresponding results on the
finitary version  of the FH from \cite{s12}, so we concentrate only
on the new moments and try to be concise whenever the material is a
straightforward extension of the previous one. We will see that the
FH of $k$-partitions is in a sense an ``iterated version'' of the DH
of $k$-partitions which is far from obvious for the the particular
case of the Wadge hierarchy, under the classical definition.

As we explain below, the FH of sets in the Baire space does coincide
with the Wadge hierarchy. To our knowledge, the extension of this
hierarchy to non-zero-dimensional spaces is new here (so far such an
extension was known only for the finitatry version of the FH
\cite{s08a,s12}). Interestingly, in our approach here the definition
of the FH of sets is in fact not simpler than for the
$k$-partitions.

Since even the definition of the FH is technically very involved, we
concentrate here on a slightly easier case of sets and
$k$-partitions of finite Borel rank and provide only proof hints for
some long proofs, appealing to the analogy with the finitary version
in \cite{s12}.

\subsection{More on the h-preorders}\label{hommore}

Here we  extend  some notions and facts from Subsection \ref{hpre}
about the $h$-preorders, in order to describe notation systems for
levels of the FHs of $k$-partitions. We omit the proofs which are
easy variations of their finitary versions in \cite{s12}.

Let  $(Q;\leq)$ be a poset. A {\em $Q$-poset}  is a triple
$(P,\leq,c)$ consisting of a countable nonempty poset $(P;\leq)$,
$P\subseteq\omega$ without infinite chains, and a labeling
$c:P\rightarrow Q$. By default, we denote the labeling in a
$Q$-poset by $c$. A {\em morphism}
$f:(P,\leq,c)\rightarrow(P^\prime,\leq^\prime,c^\prime)$ of
$Q$-posets is a monotone function
$f:(P;\leq)\rightarrow(P^\prime;\leq^\prime)$ satisfying $\forall
x\in P(c(x))\leq c^\prime(f(x)))$. Let $\widetilde{{\mathcal P}}_Q$,
$\widetilde{{\mathcal F}}_Q$ and $\widetilde{{\mathcal T}}_Q$
denote the sets of all countable $Q$-posets, $Q$-forests and
$Q$-trees, respectively.

The  {\em $h$-preorder} $\leq_h$ on ${\mathcal P}_Q$ is defined as
follows: $P\leq_h P^\prime$, if there is a morphism from $P$ to
$P^\prime$. The  quotient-posets of $\widetilde{{\mathcal P}}_Q$,
$\widetilde{{\mathcal F}}_Q$, $\widetilde{{\mathcal T}}_Q$ are
denoted $\widetilde{{\mathbb P}}_Q$, $\widetilde{{\mathbb F}}_Q$,
$\widetilde{{\mathbb T}}_Q$, respectively.  Note that for the
particular case $Q=\bar{k}$ of the antichain with $k$ elements we
obtain the preorders $\widetilde{{\mathcal P}}_k$,
$\widetilde{{\mathcal F}}_k$, $\widetilde{{\mathcal T}}_k$ from the
previous section.

Next we formulate some lattice-theoretic properties of the
$h$-preorders. By a {\em partial lower semilattice} we mean a poset
in which any two elements that have a lower bound  have a (unique)
greatest lower bound. For any poset $Q$, let $Q^\prime$ be the poset
obtained from $Q$ by adjoining the new element $\bot$ which is
smaller than all elements in $Q$.  Define the function
$s:Q\rightarrow{\mathcal P}_Q$ as follows:  $s(q)$ is the singleton
tree labeled by $q\in Q$. If $Q$ is an upper semilattice, define the
function $l:{\mathcal P}_Q\rightarrow Q$ by $l(P)=\bigcup\{c(p)\mid
p\in P\}$ where $\cup$ is the join operation in $Q$. For $Q$-posets
$P$ and $R$ (resp. $P_0,P_1,\ldots$), let $P\sqcup R$ (resp.
$\bigcup_iP_i$) denote their disjoint union (resp. countable
disjoint union), $P\sqcup R,\bigcup_iP_i\in\widetilde{{\mathcal
P}}_Q$.

For posets $P$  and $Q$ we write $P\subseteq Q$ (resp. $P\sqsubseteq
Q$) if $P$ is a substructure of $Q$ (resp. $P$ is an initial segment
of $Q$). The next assertion is essentially contained in \cite{s12}.

\begin{proposition}\label{lat}
\begin{enumerate}\itemsep-1mm

 \item If $Q$ is a partial lower semilattice then $Q^\prime$ is a lower semilattice.
 \item $(\widetilde{{\mathcal P}}_Q;\leq_h,\sqcup)$ is a distributive upper
$\sigma$-semilattice that contains $({\mathcal F}_Q;\leq_h,\sqcup)$
as a distributive upper $\sigma$-subsemilattice.
 \item If $Q$ is a partial lower semilattice then
$(\widetilde{{\mathcal P}}^\prime_Q;\leq_h)$ and $(\widetilde{{\mathcal
F}}^\prime_Q;\leq_h)$ are distributive lattices.
 \item If $Q$ is directed (i.e. any two elements have an upper bound)
then $s(Q)$ is a cofinal subset of $\widetilde{{\mathcal P}}_Q$ (i.e. any
$x\in\widetilde{{\mathcal P}}_Q$ is below $s(q)$ for some $q\in Q$).
 \item The mapping $s$ is an isomorphic embedding of $Q$ into $\widetilde{{\mathcal
 P}}_Q$.
 \item If $Q$ is a partial lower semilattice then $s$ preserves
 the greatest lower bound operations in $Q$ and $\widetilde{{\mathcal P}}_Q$
(and similarly for $\widetilde{{\mathcal F}}_Q$).
 \item If $Q$ is an upper semilattice then $l:\widetilde{{\mathcal P}}_Q\rightarrow Q$
is a homomorphism of upper semilattices and $q=l(s(q))$ for each
$q\in Q$.
 \item If $Q$ is a bqo then $(\widetilde{{\mathcal F}}_Q;\leq_h)$, $(\widetilde{{\mathcal
T}}_Q;\leq_h)$  are bqo's.
 \item Let $P,Q$ be arbitrary posets. Then $P\subseteq Q$ implies
$\widetilde{\mathbb{F}}_P\subseteq \widetilde{\mathbb{F}}_Q$, and $P\sqsubseteq Q$ implies
$\widetilde{\mathbb{F}}_P\sqsubseteq \widetilde{\mathbb{F}}_Q$.
 \end{enumerate}
 \end{proposition}

Now we can iterate the construction
$Q\mapsto\widetilde{\mathbb{F}}_Q$ starting with the antichain
$\overline{k}$ of $k$ elements $\{0,\ldots,k-1\}$. Define the
sequence $\{\widetilde{\mathcal{F}}_k(n)\}_{n<\omega}$ of preorders
by induction on $n$ as follows:
$\widetilde{\mathcal{F}}_k(0)=\overline{k}$  and
$\widetilde{\mathcal{F}}_k(n+1)=\widetilde{\mathcal{F}}_{\widetilde{\mathcal{F}}_k(n)}$.
Identifying the elements $i<k$ of $\overline{k}$ with the
corresponding minimal elements $s(i)$ of
$\widetilde{\mathcal{F}}_k(1)$, we may think that
$\widetilde{\mathcal{F}}_k(0)\sqsubseteq\widetilde{\mathcal{F}}_k(1)$.
By items (8,9) of Proposition \ref{lat},
$\widetilde{\mathcal{F}}_k(n)\sqsubseteq\widetilde{\mathcal{F}}_k(n+1)$
for each $n<\omega$, and
$\widetilde{\mathcal{F}}_k(\omega)=\bigcup_{n<\omega}\widetilde{\mathcal{F}}_k(n)$
is a wqo. For any $n\leq\omega$, let $\widetilde{\mathbb{F}}_k(n)$
be the quotient-poset of $\widetilde{\mathcal{F}}_k(n)$.

Of course, similar  constructions can be done with
$\widetilde{\mathcal{T}}$  in place of $\widetilde{\mathcal{F}}$.
The preorders $\widetilde{\mathcal{F}}_k(\omega)$,
$\widetilde{\mathcal{T}}_k(\omega)$ and the set
$\widetilde{\mathcal{T}}^\sqcup_k(\omega)$ of countable joins of
elements in $\widetilde{\mathcal{T}}_k(\omega)$, play an important
role in the study of the FH of $k$-partitions because they provide
convenient naming systems for the levels of this hierarchy (similar
to the previous section where  $\widetilde{\mathcal{F}}_k$ and
$\widetilde{\mathcal{T}}_k$ where used to name the levels of the DH
of $k$-partitions). Note that
$\widetilde{\mathcal{F}}_k(1)=\widetilde{\mathcal{F}}_k$ and
$\widetilde{\mathcal{T}}_k(1)=\widetilde{\mathcal{T}}_k$.

By Proposition \ref{lat}, for any $n<\omega$ there is an embedding
$s=s_n$ of $\widetilde{\mathbb{F}}_k(n)$ into $\widetilde{\mathbb{F}}_k(n+1)$, and
$s_{n+1}$ coincides with $s_n$ on $\widetilde{\mathbb{F}}_k(n)$. This induces
the embedding $s=\bigcup_{n<\omega}s_n$ of $\mathbb{F}_k(\omega)$
into itself such that $s$ coincides with $s_n$ on
$\widetilde{\mathbb{F}}_k(n)$ for each $n<\omega$. Similarly, for any
$n<\omega$ we have the function $l$ from $\widetilde{\mathbb{F}}_k(n+2)$ onto
$\widetilde{\mathbb{F}}_k(n+1)$ which induces the function (denoted also by
$l$) from $\widetilde{\mathbb{F}}_k(\omega)$ onto $\widetilde{\mathbb{F}}_k(\omega)$.

Define the binary operation $\ast$ on
$\widetilde{\mathcal{P}}_k(\omega)$  as follows: $F\ast G$ is the
labeled poset obtained from $G$ by adjoining a new largest (root)
element and assigning the label $F$ to that element. It is easy to
see that the operation $\ast$ respects the $h$-equivalence relation
and hence induces the binary operation on
$\widetilde{\mathbb{P}}_k(\omega)$ also denoted by $\ast$. Note that
for $F=s(i)=i<k$ we have $F\ast G=p_i(G)$ and that
$\widetilde{\mathcal{F}}_k(\omega)$ is closed under $\ast$.

We  formulate some properties of the introduced objects illustrating
a rich algebraic structure of $\widetilde{\mathcal{F}}_k(\omega)$.
They are again almost the same as their corresponding finitary
versions in \cite{s12}.

\begin{proposition}\label{dcq}
 \begin{enumerate}\itemsep-1mm
 \item For any $n$ with $1\leq n\leq\omega$, $\widetilde{\mathbb{F}}^\prime_k(n)$
is a well distributive lattice.
 \item  Any element of $\widetilde{\mathbb{F}}^\prime_k(\omega)$
is the value of a  variable-free term (countable joins are allowed) of signature
$\{\sqcup,\ast,\bot,0,\ldots,k-1\}$.
 \item For any  $0<n\leq\omega$, $\widetilde{\mathbb{T}}_k(n)$ generates
$\widetilde{\mathbb{F}}^\prime_k(n)$ under  $\sqcap$.
 \item The set of $\sigma$-join-irreducible elements of $\widetilde{\mathbb{T}}^\sqcup_k(\omega)$
 coincides with $\widetilde{\mathbb{T}}_k(\omega)$.
 \item The set of join-irreducible elements $\mathbb{T}^\sqcup_k(\omega)$ coincides with
 $\widetilde{\mathbb{T}}_k(\omega)\cup S$, where $S$ is the set of supremums
 of infinite increasing sequences in $\widetilde{\mathbb{T}}_k(\omega)$.
 \end{enumerate}
  \end{proposition}

Next we provide a characterization of
$\widetilde{\mathbb{P}}_k(\omega)$ (and the related substructures)
which is sometimes more convenient when dealing with the FH in the
next subsection. For $P,Q\in\widetilde{\mathcal{P}}_k(n)$, an {\em
explicit morphism} $\varphi:P\rightarrow Q$ is a sequence
$(\varphi_0,\ldots,\varphi_{n-1})$ where $\varphi_0$ is a morphism
from $P$ to $Q$, $\varphi_1=\{\varphi_{1,p_0}\}_{p_0\in P}$ is a
family of morphisms from $c(p_0)$ to $c(\varphi_0(p_0))$,
$\varphi_2=\{\varphi_{2,p_0,p_1}\}_{p_0\in P,p_1\in c(p_0)}$ is a
family of morphisms from $c(p_1)$ to $c(\varphi_{1,p_0}(p_1))$ and
so on (this notion makes use of the convention that $i=s(i)$ for
each $i<k$). Note that $P\leq_hQ$ iff there is an explicit morphism
from $P$ ro $Q$, and that for $n=1$ the explicit morphisms
essentially coincide with the morphisms. In the next assertion we
treat $\widetilde{\mathbb{P}}_k(\omega)$ as  the category whose
morphisms are the explicit morphisms.

For  each positive $n<\omega$, a {\em $k$-labeled $n$-preorder}
\cite{s08a} is a countable structure
$(S;d,\leq_0,\cdots,\leq_{n-1})$ where $d:S\to k$ is a $k$-partition
of $S$ and $\leq_0,\cdots,\leq_n$ are preorders on $S$ such that
$\leq_n$ is a partial order, $x\leq_{i+1}y$ implies $x\equiv_{i}y$,
the quotient-poset of $(S;\leq_0)$ has no infinite chains, for each
$x\in S$ the quotient-poset of  $([x]_0;\leq_1)$ has no infinite
chains for each $y\in [x]_0$ the quotient-poset of  $([y]_1;\leq_2)$
has no infinite chains, and so on. Let $\mathcal{S}_n(k)$ be the
category of $k$-labeled $n$-posets as objects where the morphisms
are functions that preserve the labelings and are monotone w.r.t.
all the preorders.

\begin{proposition}\label{exmor}
For any positive $n<\omega$, the categories
$\widetilde{\mathcal{P}}_k(n)$ and $\mathcal{S}_n(k)$ are
equivalent.
\end{proposition}

{\bf Proof Hint.} The equivalence is witnessed (cf. Proposition 8.8
of \cite{s12}) by the functors defined as follows. Relate to any
object $(P;\leq,c)$ of $\widetilde{\mathcal{P}}_k(n)$ the object
$P^\circ=(X;\leq_0,\ldots,\leq_{n-1},d)$ of $\mathcal{S}_n(k)$ where
$X$ is formed by the elements $p=(p_0,\ldots,p_{n-1})$ such that
$p_0\in P,p_1\in c(p_0),\ldots$, the preorders between such $p$ and
$r=(r_0,\ldots,r_{n-1})$ are defined by $p\leq_0r$ iff $p_0\leq
r_0$, $p\leq_1r$ iff $p_0=r_0$ and $p_1\leq r_1$ and so on, and the
labeling $d:X\rightarrow k$ is defined by $d(p)=c(p_{n-1})$. Relate
to any explicit morphism
$\varphi=(\varphi_0,\ldots,\varphi_{n-1}):P\rightarrow Q$ of
$\widetilde{\mathcal{P}}_k(n)$ the morphism
$\varphi^\circ:P^\circ\rightarrow Q^\circ$ by
$\varphi^\circ(p_0.p_1,\ldots)=(\varphi_0(p_0),\varphi_{1,p_0}(p_1),\ldots)$.

Conversely, relate to any object $X$ of $\mathcal{S}_n(k)$ as above
the object $X^+=(P;\leq,c)$ of $\widetilde{\mathcal{P}}_k(n)$ where
$(P;\leq)$ is the quotient-poset of $(X;\leq_0)$ and
$c([x]_0)=([x]_0;\leq_1,\ldots,\leq_{n-1},d|_{[x_0]})^+$; we can
suppose by induction that $c([x]_0)$ is an object of
$\widetilde{\mathcal{P}}_k(n-1)$ if $n>1$. Relate to any  morphism
$\psi:X\rightarrow Y$ of $\mathcal{S}_n(k)$ the explicit morphism
$\psi^+=(\varphi_0,\ldots,\varphi_{n-1}):X^+\rightarrow Y^+$ where
$\varphi_0([x]_0)=[\psi(x)]_0$,
$\varphi_{1,[x]_0}([z]_1)=[\psi(z)]_1$ for each $z\in [x]_0$ and so
on. \qed

Note that the notion of an explicit morphism does not in fact depend
on the number $n$ because the  explicit morphism
$\varphi:P\rightarrow Q$ of $P,Q\in\widetilde{\mathcal{P}}_k(n)$ is
uniquely extended to an explicit morphism of $P,Q$ considered as
objects of $\widetilde{\mathcal{P}}_k(n+1)$. Thus, we can consider
the category $\widetilde{\mathcal{P}}_k(\omega)$ with the explicit
morphism. Similarly, we can consider the category
$\mathcal{S}_\omega(k)=\bigcup_n\mathcal{S}_n(k)$ of finite because
$\mathcal{S}_n(k)$ may be considered as a subcategory of
$\mathcal{S}_{n+1}(k)$ (just add the equality relation as $\leq_n$).
In this way, we obtain the equivalence of categories
$\widetilde{\mathcal{P}}_k(\omega)$ and $\mathcal{S}_\omega(k)$.

The full subcategories $\widetilde{\mathcal{F}}_k(\omega)$ and
$\widetilde{\mathcal{T}}_k(\omega)$ of
$\widetilde{\mathcal{P}}_k(\omega)$ are then equivalent to suitable
full subcategories $\mathcal{U}_\omega(k)$ and
$\mathcal{V}_\omega(k)$ of $\mathcal{S}_\omega(k)$ (e.g., the
objects of $\mathcal{U}_\omega(k)$ are $(X;\leq_0.\leq_1\ldots,d)$
where $(X;\leq_0)$ is a forest, $([x]_0;\leq_1)$ is a forest for
each $x\in X$, and so on).

We conclude this subsection with extending Proposition \ref{t2} to
the structure $\widetilde{\mathcal{T}}^\sqcup_k(\omega)$. For this
we need the ordinal
$\xi=sup\{\omega_1,\omega_1^{\omega_1},\omega_1^{(\omega_1^{\omega_1})},\ldots\}$.
According to the Cantor normal form, any non-zero ordinal
$\alpha<\xi$ is uniquely representable in the form
$\alpha=\omega_1^{\gamma_0}\cdot
\alpha_0+\cdots+\omega_1^{\gamma_l}\cdot \alpha_l$ for some
$l<\omega$, $\alpha>\gamma_0>\cdots>\gamma_l$ and non-zero ordinals
$\alpha_0,\ldots,\alpha_l<\omega$). For
$F\in\widetilde{\mathcal{F}}_2(\omega)$, let
$\bar{F}\in\widetilde{\mathcal{F}}_2(\omega)$ be obtained from $F$
by interchanging $\{0,1\}$ in all the labels.

\begin{definition}\label{mintree}
 We define the sequence $\{T_{\alpha}\}_{\alpha<\xi}$ of
 trees in $\widetilde{\mathcal{T}}_2(\omega)$
 by induction on $\alpha$ as follows:%\\[-8mm]
 \begin{enumerate}\itemsep-1mm
 \item For $\alpha<\omega_1$, use the definition from the end of Subsection \ref{hpre}.
 \item For any non-zero ordinal $\gamma<\xi$,
 $T_{\omega_1^\gamma}=s(T_\gamma)$.
 \item For any limit uncountable ordinal $\lambda<\xi$ of countable cofinality,
 $T_\lambda=\bigsqcup_nT_{\alpha_n}$ where $\alpha_0<\alpha_1<\cdots$
 and $\lambda=sup\{\alpha_0,\alpha_1,\ldots\}$.
 \item For any  ordinal $\beta\geq\omega_1$,
$T_{\beta+1}=0\ast(T_\beta\cup \bar{T}_\beta$.
 \item For all $\delta,\gamma$ such that $1\leq \delta<\omega_1$ and $1\leq\gamma<
\xi$, $T_{\omega_1^\gamma(\delta+1)}
=T_\gamma\ast\bar{T}_{\omega_1^\gamma \delta}$.
 \item For all $\beta,\gamma$ such that $1\leq\gamma<
\xi$ and $\beta =\omega_1^{\gamma_1}\cdot\beta_1$ for some
$\beta_1>0$ and $\gamma_1>\gamma$, $T_{\beta+\omega_1^\gamma}
=T_\gamma\ast(T_\beta\sqcup\bar{T}_\beta)$.
 \item For all $\delta,\beta,\gamma$ such that $1\leq \delta<\omega_1$, $1\leq\gamma<
\xi$ and $\beta =\omega^{\gamma_1}\cdot\beta_1$ for some $\beta_1>0$
and $\gamma_1>\gamma$, $T_{\beta+\omega_1^\gamma(\delta+1)}
=T_\gamma\ast\bar{T}_{\beta+\omega_1^\gamma \delta}$.
  \end{enumerate}
 \end{definition}

The next assertion  is the extension of Proposition \ref{t2}. We
omit the proof which is very similar to that of Proposition 8.29 in
\cite{s12} which is a finitary version of this assertion.

 \begin{proposition}\label{mintree1}
 \begin{enumerate}\itemsep-1mm
 \item  For all $\alpha<\beta<\xi$, $T_\alpha$ is correctly defined up to
$\equiv_h$ (in particular, it does not depend on the choice of the
ordinals $\alpha_n$
 in Definition \ref{mintree}(3)),  $T_\alpha$ and
$\bar{T}_\alpha$ are $\leq_h$-incomparable elements of
$\widetilde{\mathcal{T}}_2(\omega)$
 that satisfy $T_\alpha\sqcup\bar{T}_\alpha<_hT_\beta$.
 \item Any $T\in\widetilde{\mathcal{T}}_2(\omega)$ is $h$-equivalent to precisely one of
 $T_\alpha,\bar{T}_\alpha$ for some $\alpha<\xi$.
 \item Any $T\in\widetilde{\mathcal{T}}^\sqcup_2(\omega)$ is $h$-equivalent to precisely one of
$T_\alpha,\overline{T}_\alpha,T_\alpha
\sqcup\overline{T}_\alpha,\bigsqcup_{\alpha<\lambda}T_\alpha$ where
$\lambda$ is a limit ordinal of countable cofinality.
 \item $T_\alpha\sqcup\overline{T}_\alpha\equiv_hT_{\alpha+1}\sqcap\overline{T}_{\alpha+1}$
and
$\bigsqcup_{\alpha<\lambda}T_\alpha\equiv_hT_{\lambda}\sqcap\overline{T}_{\lambda}$
for any limit ordinal $\lambda$ of countable cofinality, where
$\sqcap$ is a binary operation of $k$-forests inducing the infimum
operation in $\widetilde{\mathbb F}_k^\prime(\omega)$.
 \item The ranks of
$\widetilde{\mathbb{T}}_2(\omega)$ and of
$\widetilde{\mathbb{T}}^\sqcup_2(\omega)$ coincide with $\xi$.
 \end{enumerate}
  \end{proposition}

\subsection{Fine hierarchies of k-partitions}\label{fhpart}

Here we define the TH of $k$-partitions and formulate its basic
properties. The proofs are almost the same as those for the finitary
case in \cite{s12}.

Let  $X$ be a space.  By an {\em $\omega$-base} (cf.
\cite{s08a,s12}) in $X$ we mean a sequence $\mathcal {L}=\{\mathcal
{L}_n\}_{n<\omega}$ of bases such that $\mathcal
{L}_n\cup\check{\mathcal L}_n\subseteq\mathcal {L}_{n+1}$ for each
$n<\omega$. The main example of an $\omega$-base is of course
$\{\mathbf{\Sigma}^0_{n+1}\}_{n<\omega}$ but also  other examples
are interesting, in particular the finite shifts
$\{\mathbf{\Sigma}^0_{m+n}\}_{ n<\omega}$ for any fixed $1\leq
m<\omega$.

Let $P\in\widetilde{\mathcal{F}}_k(n)$ for
some positive $n<\omega$. By a {\em $P$-family over ${\mathcal
L}$} we mean a family
$\{B_{p_0},B_{p_0,p_1},\ldots,B_{p_0,\ldots,p_{n-1}}\}$ where
$p=(p_0,\dots,p_{n-1})\in P^\circ$ (see the end of Subsection
\ref{hommore}), $B_{p_0}\in{\mathcal L}_0,B_{p_0,p_1}\in{\mathcal
L}_1,\ldots,B_{p_0\ldots,p_{n-1}}\in{\mathcal L}_{n-1}$, and the sets
 $$\tilde{B}_{p_0}=B_{p_0}\setminus\bigcup\{B_r\mid r\leq_0 p\},\;
\tilde{B}_{p_0,p_1}=B_{p_0,p_1}\setminus\bigcup\{B_r\mid r\leq_1
p\},\ldots$$
 satisfy
 $$\tilde{B}_{p_0}=\bigcup\{B_{p_0,p_1}\mid p_1\in
c(p_0)\},\;\tilde{B}_{p_0,p_1}=\bigcup\{B_{p_0,p_1,p_2}\mid p_1\in
c(p_0),p_2\in c(p_1)\},\ldots.$$
 To simplify notation, we often denote families just
by $\{B_p\}$. Note that $d(p)=c(p_{n-1})$ is always in
$\widetilde{\mathcal{F}}_k(0)$=\{0,\ldots,k-1\},
 $$\tilde{B}_{p_0}=\bigcup\{\tilde{B}_{p_0,p_1}\mid p_1\in
c(p_0)\},\;\tilde{B}_{p_0,p_1}=\bigcup\{\tilde{B}_{p_0,p_1,p_2}\mid
p_1\in c(p_0),p_2\in c(p_1)\},\ldots$$
 and that for $n=1$ the
$P$-families over ${\mathcal L}$ essentially coincide with the
$P$-families of ${\mathcal L}_0$-sets in Subsection \ref{kpart}.
Obviously, $\bigcup_{p_0}B_{p_0}=\bigcup_{p\in P^\circ}\tilde{B}_p$.
We call a $P$-family $\{B_p\}$  over ${\mathcal L}$ {\em consistent}
if $d(p)=d(q)$ whenever the components $\tilde{B}_p$ and
$\tilde{B}_q$ have a nonempty intersection. Any such consistent
$P$-family determines the $k$-partition
$A:\bigcup_{p_0}B_{p_0}\rightarrow k$ where $A(x)=d(p)$ for some
(equivalently, for any)  $p\in P^\circ$ with $x\in \tilde{B}_p$; we
also say in this case that $A$ is defined by $\{B_p\}$. Note that
this $k$-partition is determined by the defining $P$-family and it
does not  depend on the number $n$ with
$P\in\widetilde{\mathcal{F}}_k(n)$.

Let ${\mathcal L}^Y(P)$ be the set of $k$-partitions
$A:Y\rightarrow k$ defined by some $P$-family over $\mathcal{L}$.
In case $Y=X$ we omit the superscript $X$ and call (temporarily)
the family $\{{\mathcal L}(P)\}_{P\in\widetilde{\mathcal{F}}_k(\omega)}$  the
FH of $k$-partitions over $\mathcal{L}$.

For $F\in\widetilde{\mathcal{F}}_k(\omega)$, by a
{\em reduced $F$-family over $\mathcal{L}$} we mean a monotone
$F$-family $\{B_p\}$ over $\mathcal{L}$ such that $B_{p_0}\cap
B_{q_0}=\emptyset$ for all incomparable $p_0,q_0\in F$,
$B_{p_0,p_1}\cap B_{p_0,q_1}=\emptyset$ for all incomparable
$p_1,q_1\in c(p_0)$ and so on. Let $\mathcal{L}^Y_r(F)$ be the set
of partial $k$-partitions defined by the reduced $F$-families
$\{B_p\}$ over $\mathcal{L}$ such that $\bigcup_{p_0}B_{p_0}=Y$.

The next assertion is a a straightforward extension of Proposition
\ref{treedh1} proved similarly to its finitary version in
\cite{s12}.

\begin{proposition}\label{treefh1}
 \begin{enumerate}\itemsep-1mm
 \item  If $A\in{\mathcal L}^Y(P)$ then $A|_Z\in{\mathcal L}^Z(P)$
for each $Z\subseteq Y$, $Z\in\mathcal{L}_0$.
 \item  Any $A\in{\mathcal L}^Y(P)$ is defined by a monotone
 $P$-family $\{C_p\}$ (monotonicity means that $C_{q_0}\subseteq C_{p_0}$ for
 $q_0\leq p_0$, $C_{p_0,q_1}\subseteq C_{p_0,p_1}$ for
 $q_1\leq p_1$ and so on).
 \item Let $f:X_1\rightarrow X$ be a morphism of $\omega$-spaces
 and $A\in{\mathcal L}^Y(P)$ in $X$.  Then
$f^{-1}(A)\in{\mathcal L}^{f^{-1}(Y)}(P)$.
 \item If $P\leq_hQ$ then ${\mathcal
L}^Y(P)\subseteq{\mathcal L}^Y(Q)$.
 \item The collection $\{{\mathcal L}(P)\mid P\in\widetilde{\mathcal{F}}_k(\omega)\}$
is well partially ordered by inclusion.
 \item  For all
$F,G\in\widetilde{\mathcal{F}}^\prime_k(\omega)$, $\mathcal{L}(F)
\cap\mathcal{L}(G)=\mathcal{L}(F\sqcap G)$.
 \item Any level ${\mathcal L}(F)$,
 $F\in\widetilde{\mathcal{F}}_k(\omega)$,
of the FH is the intersection of finitely many ``non-self-dual''
levels ${\mathcal L}(T)$, $T\in\widetilde{\mathcal{T}}_k(\omega)$.
 \item Let   $\mathcal{L}_n$ have the $\omega$-reduction property for each $n<\omega$,
$Y\in\mathcal{L}_0$ and $F\in\widetilde{\mathcal{F}}_k(\omega)$. Then
$\mathcal{L}^Y(F)=\mathcal{L}^Y_r(F)$, in particular, $\mathcal{L}(F)=\mathcal{L}_r(F)$.
 \end{enumerate}
  \end{proposition}

The next result is an infinitary version of Proposition 8.19 in
\cite{s12}. It extends the Hausdorff-Kuratowski theorem to all
levels of the FH. We call the $\omega$-base  $\mathcal{L}$ is {\em
interpolable} if, for each $n\geq1$, $\mathcal{L}_n$ has the
$\omega$-reduction property  and   the Hausdorff-Kuratowski theorem
holds for any non-zero level of $\mathcal{L}$. In particular, the
base $\mathcal{L}=\{\mathbf{\Sigma}^0_{n+1}(X)\}_{n<\omega}$ is
interpolable in any quasi-Polish space.

\begin{theorem}\label{fhpart3}
If the $\omega$-base  $\mathcal{L}$ is interpolable then any
non-empty constituent of the FH
$\{\mathcal{L}(x)\}_{x\in\widetilde{\mathbb{F}}_k(\omega)}$ is a constituent
of the FH $\{\mathcal{L}(x)\}_{x\in\widetilde{\mathbb{T}}^\sqcup_k(\omega)}$.
  \end{theorem}

Once we have the (hopefully) right definition of the FH of
$k$-partitions, some properties of the DH are extended in a
straightforward way. In particular, this applies to Propositions
\ref{dhpart3}, \ref{dhpart4} and Theorem \ref{dhpart5}. We give the
corresponding formulations.

\begin{proposition}\label{dhpart31}
Let  $\mathcal{L}=\{\mathbf{\Sigma}^0_{n+1}\}_{n<\omega}$ and
$F\in\widetilde{\mathcal{F}}_k$. Then $\{\mathcal{L}(X,F)\}_X$ is a
family of $k$-partition classes.
\end{proposition}

\begin{proposition}\label{dhpart41}
Let $X$ be a cb$_0$-space,
$\mathcal{L}=\{\mathbf{\Sigma}^0_{n+2}(X)\}_{n<\omega}$, and
$T\in\widetilde{\mathcal{T}}_k$. Then $\mathcal{L}(X,T)$ has a
principal $\mathcal{L}(T)$-TR. If $X$ is in addition
zero-dimensional and
$\mathcal{M}=\{\mathbf{\Sigma}^0_{n+1}(X)\}_{n<\omega}$ then
$\mathcal{M}(X,T)$ has a principal $\mathcal{M}(T)$-TR.
\end{proposition}

The next result extends (with essentially the same proof) Theorem
\ref{dhpart5} to the FH of $k$-partitions.

\begin{theorem}\label{dhpart51}
Let $X$ be a cb$_0$-space, $\delta:D\to X$  an admissible
representation of $X$ ($D\subseteq\mathcal{N}$), $A:X\to k$ a
$k$-partition of $X$,
$\mathcal{L}=\{\mathbf{\Sigma}^0_{n+1}\}_{n<\omega}$, and
$F\in\widetilde{\mathcal{F}}_k$. Then $A\in\mathcal{L}(X,T)$ iff
$A\circ\delta\in\mathcal{L}(D,T)$.
\end{theorem}

\begin{corollary}\label{dhpart61}
Let $X$ be a quasi-Polish space, $\delta:\mathcal{N}\to X$ an
admissible TR of $X$, $A:X\to k$ a $k$-partition of $X$,
$\mathcal{L}=\{\mathbf{\Sigma}^0_{n+1}\}_{n<\omega}$, and
$F\in\widetilde{\mathcal{F}}_k$. Then $A\in\mathcal{L}(X,T)$ iff
$A\circ\delta\in\mathcal{L}(\mathcal{N},T)$.
\end{corollary}

For the DH of $k$-partitions, it was easy to demonstrate that it
really extends the DH of sets (which coincides with the DH of
$2$-patitions by Proposition \ref{main3}). For the FH of
$2$-partitions the same task is more complicated. The reason is that
this hierarchy should generalize the Wadge hierarchy which is
defined and relatively well understood only for the Baire space (and
some other closely related spaces). Moreover, the most popular
definition of this hierarchy is in terms of the Wadge reducibility
rather than in set-theoretic terms (in fact, there is also
set-theoretic definitions \cite{wad84,lo83} but they are very
indirect and hard to deal with). Nevertheless, we claim that the FH
of 2-partitions in the Baire space coincides with the Wadge
hierarchy. We discuss this (rather informally) in the next
subsection.

%%%%%%%%%%%%%%%%%%%%%%%%%%%%%%%%%%%%%%%%%%%%%%%%%%%%%%%%%%%%%%%%%%%%%%%%
%
%
%%%%%%%%%%%%%%%%%%%%%%%%%%%%%%%%%%%%%%%%%%%%%%%%%%%%%%%%%%%%%%%%%%%%%%%%

\subsection{FH of $k$-partitions and Wadge-like reducibilities}\label{wadge}

Let us briefly recall the definition of Wadge reducibility in the
Baire space. In \cite{wad72,wad84} W. Wadge  (with a heavy use of
the Martin determinacy theorem)  proved that the structure
$({\mathbf \Delta}^1_1(\mathcal{N});\leq_W)$ is semi-well-ordered
(i.e., it is well-founded and for all $A,B\in{\mathbf B}$ we have
$A\leq_WB$ or $\overline{B}\leq_WA$. He has also computed the rank
$\nu$ of this structure which is a rather large ordinal.

In \cite{vw76,ste80} the following deep relation of the Wadge
reducibility to the separation property was established: For any
Borel set $A$ which is non-self-dual (i.e.,
$A\not\leq_W\overline{A}$) exactly one of the principal ideals
$\{X\mid X\leq_WA\}$, $\{X\mid X\leq_W\overline{A}\}$ has the
separation property.

The mentioned results give rise to the {\em Wadge hierarchy} which
is, by definition, the sequence
$\{{\bf\Sigma}_\alpha\}_{\alpha<\nu}$ of all non-self-dual principal
ideals of $({\bf \Delta}^1_1(\mathcal{N});\leq_W)$ that do not have
the separation property and satisfy for all $\alpha<\beta<\nu$ the
strict inclusion ${\bf\Sigma}_\alpha\subset{\bf\Delta}_\beta$. As
usual, we set ${\bf\Pi}_\alpha=\{\overline{A}\mid
A\in{\bf{\Sigma}_\alpha}\}$ and
${\bf\Delta}_\alpha={\bf\Sigma}_\alpha\cap{\bf\Pi}_\alpha$. Note
that the constituents of the Wadge hierarchy are precisely the
equivalence classes induced by $\leq_W$ on Borel subsets of the
Baire space (i.e., the Wadge degrees).

As shown in \cite{wad84},
${\bf\Sigma}_\alpha={\bf\Sigma}^{-1}_\alpha(\mathcal{N})$ for each
$\alpha<\omega$, i.e. the DH over open sets is an initial segment of
the Wadge hierarchy. In order to see how much finer is the Wadge
hierarchy compared with the Borel hierarchy, we mention the
equalities from \cite{wad84} relating both hierarchies:
${\mathbf\Sigma}_1={\bf\Sigma}^0_1(\mathcal{N})$,
${\mathbf\Sigma}_{\omega_1}={\bf\Sigma}^0_2(\mathcal{N})$,
${\mathbf\Sigma}_{\omega_1^{\omega_1}}={\bf\Sigma}^0_3$ and so on.
Thus, the sets  of  finite Borel rank coincide with the sets of Wadge
rank less than
$\xi=sup\{\omega_1,\omega_1^{\omega_1},\omega_1^{(\omega_1^{\omega_1})},\ldots\}$.
Note that $\xi$ is the smallest solution of the ordinal equation
$\omega_1^\varkappa=\varkappa$. Hence, we warn the reader not to
mistake ${\bf\Sigma}_\alpha$ with ${\bf\Sigma}^0_\alpha$. To give
the reader a first impression about the Wadge ordinal we note that
the rank of the preorder $({\mathbf\Delta}^0_\omega;\leq_W)$ is the
$\omega_1$-st  solution of the ordinal equation
$\omega_1^\varkappa=\varkappa$ \cite{wad84}.

The  structure of Wadge degrees is known to have the following
properties: at the bottom (i.e.,  zero) level and at the limit
levels of uncountable cofinality we have  non-self-dual pairs of
degrees; at the limit levels of countable cofinality we have
self-dual  degrees; at successor levels the self-dual degrees and
non-sel-fdual pairs alternate. Remembering Proposition
\ref{mintree1} we immediately see that the structure of Wadge
degrees of finite Borel rank is isomorphic to the structure
$\widetilde{\mathbb{T}}^\sqcup_2(\omega)$! This observation makes
more plausible our claim that the FH of $k$-partitions extends the
Wadge hierarchy.

To explain this more precisely, we note that it is possible to
relate to any $F\in\widetilde{\mathcal{T}}^\sqcup_2(\omega)$ a
$k$-partition $A_F$ of the Baire space in such a way that $A_F$ is
Wadge complete in $\mathcal{L}(F)$ where
$\mathcal{L}=\{\mathbf{\Sigma}^0_{n+1}(\mathcal{N})\}_{n<\omega}$.
Since the proof of this result is too technical for this paper, we
postpone it to a subsequent publication and only note that very
relevant particular cases are considered in \cite{s07,s07a,s11} (in
fact, from these papers only the proof for the initial segment
$\widetilde{\mathcal{T}}^\sqcup_k(2)$ is easily extracted, while for
the general result one has to employ additional jump operators from
\cite{an06,mr09} in order to relate $A_F$ to arbitrary
$F\in\widetilde{\mathcal{T}}^\sqcup_k(\omega)$. In this way one
obtains the following result showing, in particular, that the FH of
sets really extends the Wadge hierarchy of sets of finite Borel
rank:

\begin{proposition}\label{wad1}
For any $F\in\widetilde{\mathcal{T}}^\sqcup_k(\omega)$, $A_F$  is
Wadge complete in $\mathcal{L}(\mathcal{N},F)$ where
$\mathcal{L}=\{\mathbf{\Sigma}^0_{n+1}(\mathcal{N})\}_{n<\omega}$
and, moreover, $F\leq_hG$ iff
$\mathcal{L}(\mathcal{N},F)\subseteq\mathcal{L}(\mathcal{N},G)$.

In particular, for any $\alpha<\xi$ we have
$\mathcal{L}(\mathcal{N},T_\alpha)=\mathbf{\Sigma}_\alpha$ where
$T_\alpha\in\widetilde{\mathcal{T}}_2(\omega)$ is the tree  from
Proposition \ref{mintree1}.
 \end{proposition}

Our definition of the FH of $k$-partitions is thus  a natural
extension of the Wadge hierarchy to arbitrary spaces. Interestingly,
for the important particular case of quasi-Polish spaces this
hierarchy is naturally induced from the Wadge hierarchy of
$k$-partitions of the Baire space via  admissible representations.
This follows immediately from Corollary \ref{dhpart61}.

\begin{theorem}\label{wad2}
Let $X$ be a  quasi-Polish space,
$\mathcal{L}=\{\mathbf{\Sigma}^0_{n+1}(\mathcal{N})\}_{n<\omega}$,
$\mathcal{M}=\{\mathbf{\Sigma}^0_{n+1}(X)\}_{n<\omega}$, and
$F\in\widetilde{\mathcal{T}}^\sqcup_k(\omega)$. Then $A\in
\mathcal{M}(X,F)$ iff $A\circ\delta\in\mathcal{L}(\mathcal{N},F)$
where $\delta$ is some (equivalently, any) admissible TR of $X$.
 \end{theorem}

As is well known,  the structure of Wadge degrees in many natural
non-zero-dimensional cb$_0$-spaces is very complicated (see e.g.
\cite{mss12} and references therein) so it seems hopeless to
understand these structures well enough. Nevertheless, if we
slightly weaken the notion of Wadge reducibility by extending the
class of reducing functions (this is similar to the relativization
process in Computability Theory) we obtain natural versions of Wadge
reducibility which behave similar to the classical one in many
natural spaces. This also applies to reducibilities of
$k$-partition. We illustrate this with the following assertion:

\begin{proposition}\label{wad3}
Let $X$ be a quasi-Polish  space such that $dim(X)\not=\infty$. Then
the structure of $\mathbf{\Sigma}^0_3$-degrees of $k$-partitions of
$X$ of finite Borel rank is isomorphic to
$\widetilde{\mathcal{T}}^\sqcup_k(\omega)$.
 \end{proposition}

{\bf Proof Hint.}  By Proposition \ref{propgeneralhomeo}, there is a
$\mathbf{\Sigma}^0_3$-isomorphism $f$ between $X$ and $\mathcal{N}$.
Clearly, $f$ induces an isomorphism of the quotient-structures of
$(k^X;\leq_{\mathbf{\Sigma}^0_3})$ and
$(k^\mathcal{N};\leq_{\mathbf{\Sigma}^0_3})$ which preserves the
initial segments of $k$-partitions of finite Borel rank. Therefore,
it suffices to prove the assertion for $X=\mathcal{N}$. But this is
just the  $\mathbf{\Sigma}^0_3$-relativization of Proposition
\ref{wad1}. Note that similar relativizations are employed in
\cite{mr09}.
 \qed

\section{Conclusion}\label{con}

The results of this paper suggest that DST in cb$_0$-spaces (or at
least in some rich classes of cb$_0$-spaces) resembles in many
respects the classical DST in Polish spaces. Also, the methods of
classical DST seem also to work well in this context, although a
more systematic treatment of DST in cb$_0$-spaces is desirable. In
particular, the classical theory of equivalence relations and
descriptive theory of functions on cb$_0$-spaces (more complicated
than the $k$-partitions considered here) seem to be interesting.

Of course, many details on the FH of $k$-partitions in cb$_0$-spaces
should be elaborated much more carefully than in this paper. In
fact, even the structure of Wadge degrees of Borel $k$-partitions of
the Baire space should be described much more carefully; we plan to
do this in subsequent publications (of course one cannot expect to
fulfill this task in a short single paper because even much easier
case of sets is technically very involved \cite{wad84,vw76}, and the
game-theoretic technique does not work in non-zero-dimensional
spaces). The results of this paper suggest that also the FH of
$k$-partitions in arbitrary quasi-Polish spaces, and even in more
general cb$_0$-spaces, is tractable.

A special challenge is the systematic development of DST in
non-countably based spaces, in particular, in reasonable rich enough
classes of qcb$_0$-spaces. This task could require also new methods
compared with the case of cb$_0$-spaces.

Another interesting direction is the development of effective DST in
effective (in some reasonable sense) spaces. As is well known from
Computability Theory, this task is highly non-trivial even for
``simple'' spaces like the Baire space. For topologically more
complicated spaces this direction is still widely open, although it
seems of principal importance, in particular for Computable
Analysis.

%%%%%%%%%%%%%%%%%%%%%%%%%%%%%%%%%%%%%%%%%%%%%%%%%%%%%%%%%
%
%%%%%%%%%%%%%%%%%%%%%%%%%%%%%%%%%%%%%%%%%%%%%%%%%%%%%%%%%

%%%%%%%%%%%%%%%%%%%%%%%%%%%%%%%%%%%%%%%%%%%%%%%%%%%%%%%

\end{document}